\documentclass[A4paper,twoside,11pt,english]{article}
\usepackage{babel}
\usepackage{graphicx}
\usepackage{latexsym}
\usepackage{amsmath}
\usepackage{amsthm}
\usepackage{amsfonts}
\usepackage{mathrsfs}
\usepackage{geometry}
\usepackage[latin1]{inputenc}
\usepackage{amssymb}
\usepackage{makeidx}
\usepackage{xcolor}
\usepackage{makeidx}
\usepackage{tocbibind}
\usepackage{tcolorbox}
\usepackage[framemethod=TikZ]{mdframed}
\usepackage{lipsum}
\usepackage[mathscr]{eucal}
\usepackage{array}
\usepackage{comment}
\usepackage{cite}
\usepackage{hyperref}
\usepackage{algorithm2e}
\usepackage{adjustbox}
\usepackage{caption}
\usepackage{subcaption}

\title{\textbf{Image resizing by neural network operators and their convergence rate with respect to the $L^p$-norm and the dissimilarity index defined through the continuous SSIM}}
\author{\textbf{Danilo Costarelli}\thanks{Corresponding author, Paper submitted on arXiv:2501.14857v1} \quad - \quad \textbf{Mariarosaria Natale} \quad - \quad \textbf{Michele Piconi} \vspace{0.3cm}\\
Department of Mathematics and Computer Science \\
            University of Perugia\\
        1, Via Vanvitelli, 06123 Perugia, Italy    \\  
{\small {\tt danilo.costarelli@unipg.it}} \quad - \quad  {\small {\tt mariarosaria.natale@unipg.it}} \\
{\small {\tt michele.piconi@unipg.it}}}
\date{}

\providecommand{\U}[1]{\protect\rule{.1in}{.1in}}
\theoremstyle{plain}
\newtheorem{theorem}{Theorem}

\newtheorem{corollary}[theorem]{Corollary}

\newtheorem{definition}[theorem]{Definition}

\newtheorem{lemma}[theorem]{Lemma}

\newtheorem{proposition}[theorem]{Proposition}

\theoremstyle{definition}
\newtheorem{remark}[theorem]{Remark}

\newcommand{\mau}{\geq}
\newcommand{\miu}{\leq}

\newcommand{\N}{\mathbb{N}}
\newcommand{\R}{\mathbb{R}}
\newcommand{\Z}{\mathbb{Z}}
\newcommand{\disp}{\displaystyle}
\newcommand{\be}{\begin{equation}}
\newcommand{\ee}{\end{equation}}
\newcommand{\phis}{\phi_{\sigma}}

\newcommand{\xx}{{\tt{x}}}
\newcommand{\yy}{{\tt{y}}}
\newcommand{\uu}{{\tt{u}}}
\newcommand{\kk}{{\tt{k}}}

\begin{document}
\maketitle
\begin{abstract}
{\small \noindent In literature, several algorithms for imaging based on interpolation or approximation methods are available. The implementation of theoretical processes highlighted the necessity of providing theoretical frameworks for the convergence and the error estimate analysis to support the experimental setups. In this paper, we establish new techniques for deriving quantitative estimates for the order of approximation for multivariate linear operators of the pointwise-type, with respect to the $L^p$-norm and to the so-called dissimilarity index defined through the continuous SSIM. In particular, we consider a family of approximation operators known as neural network (NN) operators, that have been widely studied in the last years in view of their connection with the theory of artificial neural networks. For these operators, we first establish sharp estimates in case of $C^1$ and piecewise (everywhere defined) $C^1$-functions. Then, the case of functions modelling digital images is considered, and specific quantitative estimates are achieved, including those with respect to the mentioned dissimilarity index. Moreover, the above analysis has also been extended to $L^p$-spaces, using a new constructive technique, in which the multivariate averaged modulus of smoothness has been employed. Finally, numerical experiments of image resizing have been given to support the theoretical results.
The accuracy of the proposed algorithm has been evaluated through similarity indexes such as SSIM, likelihood index (S-index) and PSNR, and compared with other rescaling methods, including bilinear, bicubic, and upscaling-de la Vall\'ee-Poussin interpolation (u-VPI). Numerical simulations show the effectiveness of the proposed method for image processing tasks, particularly in terms of the aforementioned SSIM, and are consistent with the provided theoretical analysis.}
\end{abstract}
\medskip\noindent
{\small {\bf AMS subject classification:} 46E30, 46E36, 41A25, 41A35, 62H35, 94A08} \newline 

\noindent {\small {\bf Keywords:} Asymptotic estimates, Neural Network operators, Sobolev spaces, Lebesgue spaces, NN-algorithm, averaged modulus of smoothness, SSIM, PSNR, S-index.}
\section{Introduction}

The literature considers several algorithms for image resizing based on interpolation or approximation methods. The implementation of theoretical processes brought back the definition of several accurate and efficient numerical methods that in the last years have been strongly improved by the help of techniques of artificial intelligence and soft computing. In addition to implementation skills, theoretical knowledge is also needed to provide a convergence and error estimate analysis to support the experimental setup.

It is well-known that any static greyscale digital image can be modeled as a discontinuous, multivariate (more precisely, bivariate) signal, hence the implemented approximation or interpolation processes must be considered in their multidimensional versions, and in a general framework which includes also not necessarily continuous signals, such as $L^p$ spaces. These considerations show the necessity to derive suitable methods to study the convergence and the error of approximation of such theoretical tools, also resorting to specific (also new) metrics.

In literature, as for the convergence of multivariate approximation methods in $L^p$-setting we can find several techniques, while concerning the study of error estimates there is a lack of methods.

Indeed, the main techniques used to derive error estimates for approximation processes in the univariate setting can not be used in the case of several dimensions. More precisely, for approximation methods in case of functions of one variable one can achieve error estimates following two basic steps: 
\begin{enumerate}
    \item proving estimates for functions in Sobolev spaces $W^r_p$;
    \item extension to the whole $L^p$ using constructive density approaches, employing for instance suitable Steklov functions.
\end{enumerate} 
This procedure is widely documented, e.g., in the monograph of Sendov and Popov \cite{sendov1988averaged}. The key point of the first step is the application of a Taylor formula with integral remainder for functions belonging to $W^r_p$; it is well-known that such formula is valid in view of the equivalence (in terms of equivalent classes) between Sobolev spaces and the space $AC^{r-1}$ of functions that are absolutely continuous, with each (ordinary) derivative $f^{(i)}$ absolutely continuous, for $i=1,...,r-1$, and with $f^{(r)} \in L^p$, where  $f^{(r)}$ exists almost everywhere (see Theorem 5.4 of p. 37, in \cite{devore1993constructive}). Passing to functions of several variables, we know that the latter equivalence for the Sobolev spaces is no longer valid, hence the previous step 1) can not be replicated, making the previous strategy inapplicable. Indeed, a necessary and sufficient condition for which a function $f \in L^p(\Omega)$, $1 \miu p<+\infty$, with $\Omega \subset \R^d$, $d \mau 2$, belongs to the Sobolev spaces $W^1_p(\Omega)$ is that $f$ possesses a representative $\bar f$ that is absolutely continuous on almost all line segments in $\Omega$ parallel to the coordinate axes and whose (classical) partial derivatives belong to $L^p(\Omega)$ (see Theorem 2.14 of \cite{ziemer2012weakly}). Thus, one can provide examples of functions of several variables belonging to $W^1_p(\Omega)$ for which does not exist a continuous representative; a counterexample can be given by $f(\xx)=\|\xx\|^{-\alpha}$, defined a.e. for $\xx \in B_1 \subset \R^d$, $d-p(\alpha+1)>0$, where $B_1$ denotes the unitary ball of $\R^d$ centered at the origin, and $\|\cdot\|$ is the Euclidean norm.
\vskip0.2cm

Based on the latter considerations, the problem of deriving new strategies for proving as accurate as possible error estimates with respect to the $L^p$-norm for multivariate approximation processes arises.
\vskip0.2cm

In order to face the above issue, in this paper we propose a new implementation for image resizing based on a family of multivariate linear operators that have been deeply studied in last years, the so-called neural network (NN) operators \cite{cardaliaguet1992approximation, chen2009approximation, chen2012construction,Y21}. The name of ``NN operators'' arises from historical reasons, and it is due to the connection of the kernel function used to define such operators (generated by sigmoidal functions) with the classical theory of artificial neural networks \cite{smale2007learning,kadak2022multivariate,
kuurkova2019classification,anastassiou1997rate,
piconi2024,coroianu2024approximation,qian2022rates,
kadak2022neural}.

One of the pioneer works in this area is given by Cybenko \cite{cybenko1989approximation}, who showed the capacity of single-layer (shallow) NNs to approximate functions through the superposition of sigmoidal activation functions, using a non-constructive method based on functional analysis tools. A constructive approach to approximation was introduced through the study of (one-dimensional) NN operators, which are part of a broader class of positive linear operators. Cardaliaguet and Euvrard were among the first to propose a specific form of NN operators in \cite{cardaliaguet1992approximation}, focusing on operators with bell-shaped activation functions. While these early results provided a starting point, they had some limitations, in terms of their rate of convergence.

The latter consideration has motivated the new (more general) definition of NN operators based on sigmoidal functions given in \cite{COSP1,COSP3} in the univariate and multivariate setting, respectively. The new definition gave the possibility to overcome the main theoretical issues of the Cardaliaguet and Euvrard NN operators.
\vskip0.2cm

We recall that, the above NN operators have the following general form
\begin{equation}
\sum_{k=0}^n f(\xx_k)\, \chi_n(\xx), \quad \xx \in \Omega \subset \R^d, \tag{I}
\end{equation}
where $\chi_n$ are suitable bases functions (generated by sigmoidal functions), whose definitions will be recalled throughout the paper. Such operators are then defined using certain pointwise values of a given function $f$ at suitable nodes; it is well-known that operators of this kind are not very suitable to consider $L^p$-approximation of functions that are defined a.e., see \cite{sendov1988averaged}.

Below we describe the new strategy given in this paper to face the problem of $L^p$ estimates for the order of approximation of multivariate {\em pointwise} approximation operators.

First, we establish an asymptotic estimate for NN operators applied to both $C^1$ and piecewise $C^1$ functions (defined everywhere), with respect to the $L^p$-norm. These results are obtained under specific moment-type conditions on the generating density function and the assumption that the first-order partial derivatives are bounded by suitable positive constants. The latter condition on the first-order partial derivatives of $f$ is technical, and allows us to establish a comparison between certain line integrals (computed on segments inside the domain of the considered function) and the corresponding multiple-integral over the full domain. This relation is non-trivial, since we know that in general these kinds of inequality are not always valid.

Later, the case of image reconstruction is introduced, providing specific estimates with respect to the $L^p$-norm for piecewise constant functions, in fact providing a continuous model of (discrete) digital images. 

In this way, we obtain a general theorem which is valid under mild assumptions on the sigmoidal activation function $\sigma$; however, this generality makes the estimate not really sharp.

Two of the most common examples of functions $\sigma$ (and the corresponding density functions $\phis$) are those generated by the logistic function $\sigma_\ell$ and the ramp function $\sigma_R$.

  The first is the typical example of sigmoidal function, while the second one can be interesting in view of its relation with the so-called {\em Rectified Linear Unit function (ReLU)}, that is of central interest in the theory of artificial neural networks. 

For these very special cases, we established sharper estimates than the one achieved in the general case, especially for digital images and for the corresponding implementation for image resizing. In this context, also a pseudocode of the proposed algorithm is provided.

   In view of the above considerations on the pointwise dependence of $F_n$ on $f$, other than some difficulties to set the operators in the $L^p$ spaces, also the 
usual (integral) modulus of smoothness can not be used to evaluate their rate of convergence. For this reason, the concept of the {\em averaged modulus of smoothness} (also known as $\tau$-modulus) has been introduced in the literature as the most suitable tool that allows to estimate the order of approximation in these cases (see \cite{kolomoitsev2023sharp}).

Then, inspired by a paper by Quak \cite{quak1989multivariate}, we achieve suitable quantitative estimates for the order of approximation for the above operators for functions belonging to the whole $L^p$-spaces, using the multivariate $\tau$-modulus. The proposed proof is constructive. These estimates are more general than those established in the previous cases; however, they are less sharp.

As mentioned before, one of the main usefulness of such results is to provide a theoretical framework that justifies the performance of algorithms for image reconstruction and resizing. Numerically speaking, to concretely evaluate the accuracy of such methods, one usually resorts to discrete similarity indexes, such as the {\em Peak-Signal-to-Noise-Ratio}, the {\em Structural Similarity Index Measure} (respectively denoted by PSNR and SSIM, see, e.g., \cite{hore2010image}) or others.

Recently, Marchetti \cite{marchetti2021} introduced the so-called continuous SSIM (cSSIM) as a starting point to define a new metric for theoretically evaluate the convergence and the order of approximation of numerical methods for imaging; in \cite{marchetti_santin2022} the cases of bilinear and bicubic interpolation have been considered in details. The above-mentioned metric is the {\em continuous dissimilarity index} defined through the cSSIM, representing the continuous version of the SSIM.

Based on the previous theoretical analysis, we established the convergence and the order of approximation of the implemented algorithm for image reconstruction based on the NN operators in terms of the above mentioned continuous dissimilarity index.

Finally, at the end of the paper the introduced algorithm has been compared (through suitable numerical simulations and experiments) with some well-known interpolation methods for imaging such as  the bilinear, the bicubic, and the upscaling-de la Vall\'ee-Poussin methods (see \cite{occorsio2023image}), in order to evaluate its numerical performances.
\vskip0.2cm

Herein, we give a detailed description of the structure of the paper.

In Section \ref{s2}, we recall some basic notions on multivariate neural network operators, together with some preliminary results. Section \ref{sec2.2} is devoted to asymptotic estimates in the $L^p$ norm, for both $C^1$ and piecewise $C^1$ functions.
In Section \ref{sec2.33}, we focus on the theoretical error analysis in case of greyscale images (i.e., for piecewise constant functions), we introduce the resizing algorithm based on NN operators providing its pseudo-code.
Furthermore, we also consider in details the special cases of NN operators activated by the logistic and ramp sigmoidal activation functions, for which we obtain a sharper estimate than the one in the general case. 
Section \ref{secLp} presents convergence rates in $L^p$ spaces, obtained via a constructive density method based on the averaged modulus of smoothness (or $\tau$-modulus), allowing us to extend the estimates to the full $L^p$ spaces. Here a crucial role is played by the application of a Riesz-Thorin interpolation theorem for linear operators and $\tau$-moduli.
In Section \ref{3.1}, we first recall the SSIM evaluation metric, along with its continuous version, with the aim to establish a precise connection between image similarity, measured by the above indexes, and the $L^2$-norm. Finally, the convergence and the order of approximation of the NN operators applied in the case of digital images with respect to the continuous dissimilarity index is established.
Sections \ref{3.2} and \ref{s4} discuss alternative evaluation metrics for the numerical methods for image resizing, together with others konwn in the literature and used in this paper for comparison. 
Section \ref{s5} is devoted to numerical experiments and results.  
Finally, Section \ref{final} contains concluding remarks and possible future lines of research.

%
\section{Some recall on the multivariate neural network operators} \label{s2}

We recall that a function $\sigma: \R \to \R$ is called a {\em sigmoidal function} if and only if $\lim_{x \to -\infty} \sigma(x) = 0$ and $\lim_{x \to +\infty} \sigma(x) = 1$. In what follows, we consider non-decreasing functions $\sigma$, satisfying the following assumptions
\begin{description}
	\item[$(\Sigma 1)$] $g_{\sigma}(x) := \sigma(x) - 1/2$ is an odd function;
	\item[$(\Sigma 2)$] $\sigma \in C^2(\R)$ is concave for $x \geq 0$;
	\item[$(\Sigma 3)$] $\sigma(x) = \mathcal{O}(|x|^{-1-\alpha})$ as $x \to -\infty$, for some $\alpha > 0$.
\end{description}

\noindent For any $\sigma$  satisfying all the above assumptions, we define
\begin{equation*}
   \phi_{\sigma}(x) := \frac{1}{2}[\sigma(x+1) - \sigma(x-1)], \quad  x \in \R,
\end{equation*}
and its discrete absolute moment of order $\beta\geq 0$ given by
\begin{equation} \label{momento-d=1}
	M_{\beta}(\phi_{\sigma}):=\sup_{x\in\R}\sum_{k\in\Z}\phi_{\sigma}(x-k)|x-k|^{\beta}.
\end{equation}
For an overview of the key properties of $\phi_{\sigma}$, we refer to \cite{COSP1}. We now recall the definition of the multivariate function generated by $\sigma$
\begin{equation*}
  \Psi_{\sigma}({\tt{x}}) := \phi_{\sigma}(x_1) \cdot \phi_{\sigma}(x_2) \cdot \dots \cdot \phi_{\sigma}(x_d), 
  \quad  {\tt{x}} := (x_1, ..., x_d) \in \R^d.	
\end{equation*} 
For preliminary properties of $\Psi_{\sigma}$, the reader can see Lemma 2.3, Lemma 2.4, Lemma 2.5 and Lemma 2.7 of \cite{COSP3}.\\

From now on, we denote by $I:=\prod_{i=1}^{d} [a_i,b_i]$ the $d$-dimensional interval on $\R^d$, with $a_i,b_i\in\R$, $i=1,\cdots,d$, and by $C^0(I)$ the space of all continuous real-valued functions defined on $I$, equipped with the usual sup-norm $\| \cdot \|_{\infty}$. Further, let introduce the following set of indexes
\begin{equation*}
\mathcal{A}_n:=\{{\tt{k}}\in\Z^d \;:\; \lceil {na_i \rceil } \leq k_i \leq \lfloor {nb_i \rfloor},\; i=1,\cdots,d \}, \quad n\in\N.
\end{equation*}
\begin{lemma}[\cite{COSP3}] \label{lemma1} \rm
Let ${\tt{x}}\in I$ and $n\in\N^+$, then
\begin{equation*}
   m_{0,0}^n(\Psi_\sigma,n{\tt{x}}):= \sum_{{\tt{k}}\in \mathcal{A}_n}\Psi_\sigma(n{\tt{x}}-{\tt{k}})\geq \left[\phis(1)\right]^d>0.
\end{equation*}
\end{lemma}
Note that, in Lemma \ref{lemma1} it turns out that $\phis(1)>0$ in view of the properties of $\sigma$.
Further, we can also define the \textit{discrete absolute moment of $\Psi_{\sigma}$ of order $\beta\geq0$} (extending (\ref{momento-d=1}) to the multivariate setting) as follows
\begin{equation*}
M_{\beta}(\Psi_{\sigma}):=\sup_{{\tt{x}}\in\R^d}\sum_{{\tt{k}}\in \Z^d} \Psi_\sigma({\tt{x}}-{\tt{k}})\left\|{\tt{k}}-{\tt{x}}\right\|^{\beta},
\end{equation*}
where $\left\|\cdot\right\|$ denotes the usual Euclidean norm on $\R^d$.
It is well-known that 
\begin{equation}\label{momento}
M_0(\Psi_\sigma)=1, \quad and \quad   M_{\beta}(\Psi_{\sigma})<+\infty, \qquad 0\leq \beta<\alpha, \quad \quad d \mau 1,
\end{equation}
see, e.g., \cite{CJat23}. We now establish the following lemma, that is a slight variant of Lemma 2.5 of \cite{COSP3}.
\begin{lemma}\label{lemmanuovo}
    For every $\gamma > 0$, $0\leq\theta<1$, and $0 < \nu < \alpha$, there holds
\begin{equation*}
\sum_{\|{\tt{x}} - {\tt{k}}\| > \gamma n^{1-\theta}} \Psi_\sigma({\tt{x}} - {\tt{k}}) = O\left(n^{-\nu(1-\theta)}\right), \quad \text{as } n \to +\infty,
\end{equation*}
where $\alpha > 0$ is the constant appearing in condition $(\Sigma 3)$.
\end{lemma}
\begin{proof}
     For every $\xx \in 
\R^d$, we can write
$$
    \sum_{\kk \in \Z^d} \Psi_{\sigma}(\xx-\kk) 
  = \sum_{k_j \in \Z} \phi_{\sigma}(x_j-k_j) \left[ \, \sum_{\kk_{[j]} \in \Z^{d-1}} 
        \Psi^{[j]}_{\sigma}(\xx_{[j]}-\kk_{[j]}) \, \right],
$$
where 
$$
  \Psi^{[j]}_{\sigma}(\xx_{[j]}-\kk_{[j]}) := \phi_{\sigma}(x_1-k_1) \cdot \, 
        ... \, \cdot \phi_{\sigma}(x_{j-1}-k_{j-1}) \cdot 
    \phi_{\sigma}(x_{j+1}-k_{j+1}) \cdot \, ... \, \cdot \phi_{\sigma}(x_d-k_d),
$$
being $\xx_{[j]} := (x_1, ..., x_{j-1}, x_{j+1}, ..., x_d)$ $\in $$\R^{d-1}$, 
$\kk_{[j]} := (k_1, ..., k_{j-1}, k_{j+1}, ..., k_d) \in \Z^{d-1}$, for every $j
= 1,...,d$.  Let now $0<\nu<\alpha$. Recalling that in $\R^d$ all norms are equivalent, there exists $C>0$ such that
$$
\|\xx\|_m\, :=\ \max\{ |x_i|:\ i=1,\, ...,\, d  \}\ \leq\ C\, \| \xx \|,
$$
hence for $n \in \N^+$ sufficiently large
\begin{equation*}
    \begin{split}
        &\sum_{\| \xx - \kk \| > \gamma n^{1-\theta}} \Psi_{\sigma}(\xx-\kk)\ \leq \sum_{\| \xx - \kk \|_m > (\gamma n^{1-\theta})/C} \Psi_{\sigma}(\xx-\kk) 
\\
   &\miu \sum^d_{j=1} \, \left\{ \, \sum_{|x_j - k_j| > (\gamma n^{1-\theta})/C} 
       \phi_{\sigma}(x_j-k_j) \left[ \, \sum_{\kk_{[j]} \in \Z^{d-1}} 
     \Psi^{[j]}_{\sigma}(\xx_{[j]}-\kk_{[j]}) \, \right] \, \right\}\\
   &\miu \sum^d_{j=1} \, \left\{ \, \sum_{|x_j - k_j| > (\gamma n^{1-\theta})/C} 
     \phi_{\sigma}(x_j-k_j) \, \right\}
            \\
      &= \sum^d_{j=1} \, \left\{ \, \sum_{|x_j - k_j| > (\gamma n^{1-\theta})/C} 
  \phi_{\sigma}(x_j-k_j) \, \frac{|x_j - k_j|^{\nu}}{|x_j - k_j|^{\nu}} \right\}\\
   & < \frac{C^\nu}{\gamma^{\nu} \, n^{\nu(1-\theta)}} 
         \sum^d_{j=1} \, \left\{ \, \sum_{|x_j - k_j| > (\gamma n^{1-\theta})/C} 
         \phi_{\sigma}(x_j-k_j) \, |x_j - k_j|^{\nu} \right\}
 \end{split}
\end{equation*}
\begin{equation*}
\begin{split}
 &\leq d\frac{M_\nu(\phis)}{\gamma^{\nu}n^{\nu(1-\theta)}}<+\infty,
\end{split}
\end{equation*}
in view of (\ref{momento}) applied with $d=1$. This completes the proof.
\end{proof}
Now, we recall the definition of the classical neural network (NN) operators, that will be used in the following sections.
\begin{definition}\label{def_op} \rm
  Let $f: I \to \R$ be a given function, and $n \in \N^+$. The 
{\em multivariate} NN {\em operators} $F_n(f,\xx)$, activated 
by the sigmoidal function $\sigma$, and acting on $f$, are defined by
\begin{equation*}
	 F_n(f,\xx) := \frac{\disp \sum_{{\tt{k}}\in \mathcal{A}_n}f\left({\tt{k}}/n \right) \Psi_{\sigma}(n{\tt{x}}-{\tt{k}})}
      {\disp \sum_{{\tt{k}}\in \mathcal{A}_n}
           \Psi_{\sigma}(n{\tt{x}}-{\tt{k}})}, 
\end{equation*}
for every $\xx \in I$ and ${\tt{k}}/n := (k_1/n, ..., k_d/n)$.
\end{definition}

In view of Lemma \ref{lemma1}, the operators $F_n$ are well-defined, e.g., if the function $f$ is bounded. It is easy to see that $\left|F_n(f,{\tt{x}})\right|\leq \|f\|_\infty$. We also recall a first fundamental result, which asserts that the NN operators $F_n$ converge pointwise and uniformly in case of continuous functions on $I$, as $n \to +\infty$ (see Theorem 3.2 of \cite{COSP3}).

\section{Asymptotic estimates for multivariate operators: new approaches} \label{sec2.2}

From now on, we denote by $L^p(I)$, with $1\le p<+\infty$, the spaces of all measurable functions $f$ such that $\left|f\right|^p$ is integrable on $I$, endowed with the $L^p$-norm $\|f\|_p:=\left(\int_{I}\left|f({\tt{x}})\right|^pd{\tt{x}}\right)^{1/p}$. Moreover, we recall the definition of Sobolev spaces
\begin{equation*}
    W_p^r(I):=\left\{f\in L^p(I) \,\,\,\text{:}\,\,\, D^{\tt h}f\in L^p(I),\,\,\forall{\tt h}\in \N^d, |{\tt h}| \le r,\,\, r\in\N \right\},
\end{equation*}
where $D^{\tt h}f=\displaystyle\frac{\partial^{|\tt h|}f}{\partial x_1^{h_1}\cdots\partial x_d^{h_d}}$, $D^0f=f$, with $|{\tt h}|:=h_1+\dots+h_d$, denote, as usual, the weak partial derivatives of $f$. A norm in this space can be defined by
\[
\|f\|_{r,p}:=\sum_{0 \miu|{\tt h}|\leq r}\|D^{\tt h}f\|_p.
\]
The study of families of pointwise operators in the $L^p$-setting (such as $F_n$) is non trivial, since $L^p$-convergence in general fails. Indeed, if we consider for instance the function $h:I \to \R$
\be
h(\xx)\, :=\, \begin{cases}
1, \quad I \cap \mathbb{Q}^d,\\
0, \quad \mbox{otherwise},
\end{cases}
\ee
one has $F_n(h,\xx)=1$, $\xx \in I$, while $\|F_n(h,\cdot) - h(\cdot)\|_p=|I|^{1/p}\neq 0$, for every $n \in \N$, where $|I|$ denotes the Lebesgue measure of $I$.\\

The result below serves as a preliminary step to achieve the main aim of the present paper, that is to extend the error analysis to the whole $L^p(I)$. For the sake of simplicity, we consider only the two dimensional case $d=2$, however, the results can be easily extended to domains of higher dimension.\\

From now on, we will always consider functions $f:I\to \R$ (and its partial derivatives) as extended to the whole $\R^2$ with periodicity with respect of each variable.
\begin{theorem}\label{th_nuovo}
Let $\sigma$ be a sigmoidal function satisfying condition $(\Sigma 3)$ with $\alpha>p$, with $1\leq p< +\infty$. Moreover, we consider any function $f\in C^{1}(I)$ for which there exist positive constants $M_f$, $\beta_f>0$, such that
\begin{equation} \label{lower-bound-D}
 M_f\, \geq\, \left|\frac{\partial f}{\partial x_i}(\tt{x}) \right|^p\, \geq\, \beta_f\, >0, \quad {\tt{x}} \in I, \quad \mbox{with} \quad i=1,2.  
\end{equation}
Then there exists a positive constant $K_1$ (independent on $n$ and $f$) such that 
\begin{equation*}
    \|F_n(f,\cdot)-f(\cdot)\|_p\le \frac{K_1}{n}\, \left({M_f \over \beta_f}\right)^{1/p}\, \|f\|_{1,p}\;,\qquad n\in\N.
\end{equation*}
\end{theorem}
\begin{proof}
Let $f\in C^{1}(I)$ be fixed. The following first order Taylor formula with integral remainder holds (see \cite{folland2005higher})
\begin{equation}\label{formulataylor}
f({\tt{u}})=f({\tt{x}})+R_{1}\left({\tt{u}},{\tt{x}}\right):=f({\tt{x}})+\sum_{i=1}^2\left(u_i-x_i\right)\int_{0}^{1} \frac{\partial f}{\partial x_i}\left({\tt{x}}+t({\tt{u}}-{\tt{x}})\right)dt,
\end{equation}
${\tt{x}},{\tt{u}}\in I$. Now, using (\ref{formulataylor}), for every fixed ${\tt{x}}\in I$ and every positive integer $n$, we can write what follows
\begin{equation*}
\begin{split}
F_{n}(f,{\tt{x}})=f({\tt{x}})+ \frac{\displaystyle\sum_{{\tt{k}}\in \mathcal{A}_n}\Biggl\{ \sum_{i=1}^2\left(\frac{k_i}{n}-x_i\right)\int_{0}^{1} \frac{\partial f}{\partial x_i}\left({\tt{x}}+t\left(\frac{{\tt{k}}}{n}-{\tt{x}}\right)\right)dt\Biggr\}\Psi_\sigma(n{\tt{x}}-{\tt{k}})}{ \displaystyle\sum_{{\tt{k}}\in \mathcal{A}_n}\Psi_\sigma(n{\tt{x}}-{\tt{k}})}.
\end{split}
\end{equation*}
Thus, using the above expression for the operator $F_n$, we obtain for the fixed $1 \leq p <+\infty$
\begin{equation*}
\begin{split}
&\int_{I}|F_n (f,{\tt{x}})-f({\tt{x}})|^{p}d{\tt{x}}\leq \int_{I}\left|\frac{\displaystyle\sum_{{\tt{k}}\in \mathcal{A}_n}\Biggl\{ \sum_{i=1}^2\left(\frac{k_i}{n}-x_i\right)\int_{0}^{1} \frac{\partial f}{\partial x_i}\left({\tt{x}}+t\left(\frac{{\tt{k}}}{n}-{\tt{x}}\right)\right)dt\Biggr\}\Psi_\sigma(n{\tt{x}}-{\tt{k}})}{ \displaystyle\sum_{{\tt{k}}\in \mathcal{A}_n}\Psi_\sigma(n{\tt{x}}-{\tt{k}})}\right|^{p} d{\tt{x}}.
\end{split}
\end{equation*}
By using Jensen inequality twice, Lemma \ref{lemma1}, and again the convexity of $|\cdot|^p$, we have
\begin{equation*}
\begin{split}
&\int_{I}|F_n (f,{\tt{x}})-f({\tt{x}})|^{p}d{\tt{x}}\\
&\leq \int_{I}    \frac{\displaystyle\sum_{{\tt{k}}\in \mathcal{A}_n}\left|\sum_{i=1}^2\left(\frac{k_i}{n}-x_i\right)\int_{0}^{1} \frac{\partial f}{\partial x_i}\left({\tt{x}}+t\left(\frac{{\tt{k}}}{n}-{\tt{x}}\right)\right)dt\right|^{p}\Psi_\sigma(n{\tt{x}}-{\tt{k}})}{\displaystyle \sum_{{\tt{k}}\in \mathcal{A}_n}\Psi_\sigma(n{\tt{x}}-{\tt{k}})}\,  d{\tt{x}}
\end{split}
\end{equation*}
\begin{equation*}
\begin{split}
& \leq \frac{1}{\left[\phis(1)\right]^{2}}\int_{I} \sum_{{\tt{k}}\in \mathcal{A}_n}\left|\sum_{i=1}^2\left(\frac{k_i}{n}-x_i\right)\int_{0}^{1} \frac{\partial f}{\partial x_i}\left({\tt{x}}+t\left(\frac{{\tt{k}}}{n}-{\tt{x}}\right)\right)dt\right|^{p}\Psi_\sigma(n{\tt{x}}-{\tt{k}})d{\tt{x}}\\
&\leq \frac{2^{p-1}}{\left[\phis(1)\right]^{2}}\int_{I}\sum_{{\tt{k}}\in \mathcal{A}_n}\sum_{i=1}^2\left|\frac{k_i}{n}-x_i\right|^p\left|\int_{0}^{1} \frac{\partial f}{\partial x_i}\left({\tt{x}}+t\left(\frac{{\tt{k}}}{n}-{\tt{x}}\right)\right)dt\right|^{p}\Psi_\sigma(n{\tt{x}}-{\tt{k}})d{\tt{x}}\\
&\leq \frac{2^{p-1}}{n^p\left[\phis(1)\right]^{2}}\int_{I}\sum_{{\tt{k}}\in \mathcal{A}_n}\sum_{i=1}^2\left|k_i-n x_i\right|^p\left[\int_{0}^{1} \left|\frac{\partial f}{\partial x_i}\left({\tt{x}}+t\left(\frac{{\tt{k}}}{n}-{\tt{x}}\right)\right)\right|^pdt\right]\Psi_\sigma(n{\tt{x}}-{\tt{k}})d{\tt{x}}.
\end{split}
\end{equation*}
Being $I$ a bounded domain of $\R^2$, there exists a suitable $C>0$ sufficiently large such that
\begin{equation*}
    \Gamma_{{\tt{x}}}:=\left\{{\tt{x}}+tC(\cos\theta,\sin\theta), \, t\in[0,1],\,\theta\in[0,2\pi]\right\}
\end{equation*}
contains $I$ for every ${\tt{x}}\in I$, i.e., $I\subset \Gamma_{\tt{x}}$, for every ${\tt{x}}\in I$. Hence, we have
\begin{equation} \label{comparison-different-integrals}
\begin{split}
\int_{0}^{1}\left|\frac{\partial f}{\partial x_i}\left({\tt{x}}+t\left(\frac{{\tt{k}}}{n}-{\tt{x}}\right)\right)\right|^p dt\, \le\, {M_f \over \beta_f}\int_0^{2\pi}d\theta\int_{0}^{1}\left|\frac{\partial f}{\partial x_i}\left({\tt{x}}+tC\left(\cos\theta,\sin\theta\right)\right)\right|^p dt. 
\end{split}
\end{equation}
Indeed, it is clear that for any fixed $\tt{x} \in I$, there exists a suitable $\widetilde{\theta} \in [0,2\pi]$ (depending on $\tt{x}$ and $\tt{k}/n$) such that the segment connecting $\tt{x}$ and $\tt{k}/n$ is contained in the set
$$
\Gamma^{\widetilde{\theta}}_{{\tt{x}}}:=\left\{{\tt{x}}+tC(\cos \widetilde{\theta} ,\sin \widetilde{\theta} ), \, t\in[0,1]\right\}.
$$
This immediately implies that
$$
\int_{0}^{1}\left|\frac{\partial f}{\partial x_i}\left({\tt{x}}+t\left(\frac{{\tt{k}}}{n}-{\tt{x}}\right)\right)\right|^p dt\, \le\, \int_{0}^{1}\left|\frac{\partial f}{\partial x_i} \left({\tt{x}}+tC(\cos \widetilde{\theta} ,\sin \widetilde{\theta})\right)\right|^p dt.
$$
At the same time, for every fixed $t\in[0,1]$, by the integrability of the function $\left|\frac{\partial f}{\partial x_i}\left({\tt{x}}+tC\left(\cos\theta,\sin\theta\right)\right)\right|^p$ in the interval $[0,2\pi]$ with respect to the variable $\theta$, and by condition (\ref{lower-bound-D}), using the integral-mean theorem we know that there exists a suitable positive value $\xi=\xi({\tt{x}},{\tt{k}}/n,t)\in [\beta_f,M_f]$, such that
$$
\int_0^{2\pi} \left|\frac{\partial f}{\partial x_i}\left({\tt{x}}+tC\left(\cos\theta,\sin\theta\right)\right)\right|^p\, d\theta\ =\ 2\pi\, \xi.
$$
Thus, using the Fubini-Tonelli theorem we can get
\begin{equation*}
    \begin{split}  
&\int_{0}^{1}\left|\frac{\partial f}{\partial x_i}\left({\tt{x}}+t\left(\frac{{\tt{k}}}{n}-{\tt{x}}\right)\right)\right|^p dt\, \le\, \int_{0}^{1}\left|\frac{\partial f}{\partial x_i}\left({\tt{x}}+tC(\cos \widetilde{\theta} ,\sin \widetilde{\theta})\right)   \right|^p\cdot {\xi \over \xi}\, dt 
\\
&=\  \int_{0}^{1}\left|\frac{\partial f}{\partial x_i}\left({\tt{x}}+tC(\cos \widetilde{\theta} ,\sin \widetilde{\theta})\right) \right|^p\, {1 \over 2\pi \xi} \left[  \int_0^{2\pi} \left|\frac{\partial f}{\partial x_i}\left({\tt{x}}+tC\left(\cos\theta,\sin\theta\right)\right)\right|^p\, d\theta \right]  dt
\\
&\leq\ {M_f \over \beta_f\, 2\pi} \int_{0}^{1} \left[  \int_0^{2\pi} \left|\frac{\partial f}{\partial x_i}\left({\tt{x}}+tC\left(\cos\theta,\sin\theta\right)\right)\right|^p\, d\theta \right]\, dt
\\
&=\ {M_f \over \beta_f\, 2\pi} \int_0^{2\pi} d\theta \int_{0}^{1}  \left|\frac{\partial f}{\partial x_i}\left({\tt{x}}+tC\left(\cos\theta,\sin\theta\right)\right)\right|^p\, dt.
\end{split}
\end{equation*}
Therefore
\begin{equation*}
\begin{split}
&\int_{I}|F_n (f,{\tt{x}})-f({\tt{x}})|^{p}d{\tt{x}}
\\
&\leq {M_f \over \beta_f\, 2\pi }\frac{2^{p-1}}{n^{p}\left[\phis(1)\right]^{2}}\int_{I}\sum_{{\tt{k}}\in \mathcal{A}_n}\sum_{i=1}^2\left|k_i-nx_i\right|^p\left[\int_0^{2\pi}d\theta\int_{0}^{1}\left|\frac{\partial f}{\partial x_i}\left({\tt{x}}+tC\left(\cos\theta,\sin\theta\right)\right)\right|^p dt\right] \Psi_\sigma(n{\tt{x}}-{\tt{k}})d{\tt{x}}\\
&={M_f \over \beta_f\, 2\pi }\frac{2^{p-1}}{n^{p}\left[\phis(1)\right]^{2}}\int_{I} \sum_{i=1}^2\left\{ \left[\int_0^{2\pi}d\theta\int_{0}^{1}\left|\frac{\partial f}{\partial x_i}\left({\tt{x}}+tC\left(\cos\theta,\sin\theta\right)\right)\right|^p dt\right]\sum_{{\tt{k}}\in \mathcal{A}_n}\left|k_i-nx_i\right|^p\phi_\sigma(nx_i-k_i) \right.
\\
&\left. \qquad\times\prod_{j=1, j\neq i}^{2}\phi_\sigma(nx_j-k_j)\right\}\, d{\tt{x}}
\\
&\leq {M_f \over \beta_f\, 2\pi} \frac{ 2^{p-1}M_{p}\left(\phis\right)}{n^{p}\left[\phis(1)\right]^{2}}\int_{I}\sum_{i=1}^2\left[\int_0^{2\pi}d\theta\int_{0}^{1}\left|\frac{\partial f}{\partial x_i}\left({\tt{x}}+tC\left(\cos\theta,\sin\theta\right)\right)\right|^p dt\right] d{\tt{x}},
\end{split}
\end{equation*}
since $\prod_{j=1, j\neq i}^{2}\phi_\sigma(nx_j-k_j)\leq 1$.
Finally, by using Fubini-Tonelli theorem, and by observing that $\left\|\frac{\partial f}{\partial x_i}\left({\tt{\cdot}}\pm t C(\cos\theta, \sin \theta)\right)\right\|_p=\left\|\frac{\partial f}{\partial x_i}\right\|_p$, with $t\in [0,1]$, $\theta\in [0,2\pi]$, we obtain
\begin{equation*}
\begin{split}
\int_{I}|F_n (f,{\tt{x}})-f({\tt{x}})|^{p}d{\tt{x}}&={M_f \over \beta_f\, 2\pi}\frac{ 2^{p-1} M_{p}\left(\phi_\sigma\right)}{n^{p}\left[\phis(1)\right]^{2}}\sum_{i=1}^2\left\|\frac{\partial f}{\partial x_i}\right\|_{p}^{p}\biggl[\int_0^{2\pi}d\theta\int_0^1dt \biggr]\\
&= {M_f \over \beta_f} \frac{ 2^{p-1} M_{p}\left(\phi_\sigma\right)}{n^{p}\left[\phis(1)\right]^{2}}\sum_{i=1}^2\left\|\frac{\partial f}{\partial x_i}\right\|_{p}^{p}
\\
&=: {M_f \over \beta_f} \frac{K_1^p}{n^{p}}\sum_{i=1}^2\left\|\frac{\partial f}{\partial x_i}\right\|_{p}^{p} < +\infty,
\end{split}
\end{equation*}
where $K_1$ is finite in view of (\ref{momento}); therefore, we obtain the desired result.
\end{proof}
We remark that in (\ref{comparison-different-integrals}) we established a comparison criterion between line and multiple integrals, that in general is not valid.
\begin{remark} 
It should be noted that the {\em multidimensional} ($d>1$) Taylor formula with integral remainder (\ref{formulataylor}) is valid only for $C^1$-functions. 
As discussed in the introduction of this paper, and unlike the one-dimensional case ($d=1$), in the multidimensional setting, it is not true, in general, that for every function $f\in W^1_p(I)$ there exists an absolutely continuous function $g: I \to \R$ with $D^{{\tt{h}}}g\in L^p(I)$, $|{\tt{h}}|=1$, such that $g=f$ a.e. on $I$. This is a remarkable difference between functions of one and several variables that allows us to use (\ref{formulataylor}) for $f \in W^1_p(I)$ only when $d=1$, but not in general when $d>1$.
\end{remark}
We can prove the following.
\begin{theorem}\label{thnuovo2} 
Let $\sigma$ be a sigmoidal function satisfying condition $(\Sigma 3)$ with $\alpha>p$, $1 \miu p <+\infty$. Let $f:I\to \mathbb{R}$ be a given bounded function, with the domain $I$ that can be decomposed as
\be \label{dominio1}
I=\bigcup_{\substack{i=0,\dots,q_1-1 \\ j=0,\dots,q_2-1}}R_{ij},
\ee
where $l_0=a_1$, $t_0=a_2$, $l_{q_1}=b_1$, $t_{q_2}=b_2$, and
\begin{equation}\label{dominio2}
\begin{split}
R_{ij}&:=[l_i,l_{i+1})\times[t_j,t_{j+1}),\quad i=0,\dots,q_1-2,\quad j=0,\dots,q_2-2,\\
R_{q_1-1\ j}&:=[l_{q_1-1}, l_{q_1}]\times[t_j,t_{j+1}), \quad j=0,\dots,q_2-2,\\
R_{i\ q_2-1}&:=[l_i,l_{i+1})\times[t_{q_2-1},t_{q_2}], \quad i=0,\dots,q_1-2,\\
R_{q_1-1\ q_2-1}&:= [l_{q_1-1}, l_{q_1}]\times[t_{q_2-1},t_{q_2}],
\end{split}
\end{equation}
such that $f\in C^1(R_{ij}^\circ)$, and satisfies condition (\ref{lower-bound-D}) on every $R_{ij}$, with suitable constants $M_f$ and $\beta_f$, for every $i=0,\dots,q_1-1$, $j=0,\dots,q_2-1$. Then
\begin{equation*}
\|F_n(f,\cdot)-f(\cdot)\|_p\ \leq\ (2^{p-1}\, q_1\, q_2)^{1/p}\left( \left({M_f \over \beta_f}\right)^{1/p} \frac{K_1}{n}\left\|f\right\|_{1,p}+K_2\frac{\|f\|_\infty}{n^{\nu(1-\theta)/p}}+8^{1/p}\, \frac{\|f\|_\infty}{n^{\theta/p}} \right),
\end{equation*}
for $0<\nu<\alpha$ and $0<\theta<1$, where $K_1>0$ arises from Theorem \ref{th_nuovo}.
\end{theorem}
\begin{proof} 
First, we can observe that the operator $F_n(f,\cdot)$ is well-defined in view of the boundedness of $f$ and since $f$ is defined everywhere on $I$.  For $0<\theta<1$, we define
\[
S_{ij}:=\left[l_i+\frac{1}{n^\theta},l_{i+1}-\frac{1}{n^\theta}\right]\times \left[t_j+\frac{1}{n^\theta},t_{j+1}-\frac{1}{n^\theta}\right]\subset R_{i,j},
\]
such that every $S_{ij}^\circ\ne\emptyset$ for $n\in\N$ sufficiently large. Then for ${\tt{x}}\in R_{ij}^\circ$, with $i=0,\dots,q_1-1$, $j=0,\dots,q_2-1$,  we have
$$
\int_{R_{ij}}\left|F_n(f,\xx)-f(\xx)\right|^p\, d\xx\ \leq\ \left\{\int_{S_{ij}}+\int_{R_{ij}\setminus S_{ij}}\right\}\left|F_n(f,\xx)-f(\xx)\right|^p\, d\xx\ =:\ A + B.
$$
To estimate $A$ we consider
\begin{align}\label{spezzato}
F_n(f,\xx)-f(\xx) &= F_n(f,\xx)-f(\xx)\cdot\frac{\displaystyle\sum_{{\tt{k}}\in \mathcal{A}_n}\Psi_\sigma(n{\tt{x}}-{\tt{k}})}{\displaystyle\sum_{{\tt{k}}\in \mathcal{A}_n}\Psi_\sigma(n{\tt{x}}-{\tt{k}})}=\frac{\displaystyle\sum_{{\tt{k}}\in \mathcal{A}_n}\left(f\left(\frac{{\tt{k}}}{n}\right)-f(\xx)\right)\Psi_\sigma(n{\tt{x}}-{\tt{k}})}{\displaystyle\sum_{{\tt{k}}\in \mathcal{A}_n}\Psi_\sigma(n{\tt{x}}-{\tt{k}})}\notag\\
&=\frac{\displaystyle\sum_{{\tt{k}}/n\in R_{ij}^\circ}\left(f\left(\frac{{\tt{k}}}{n}\right)-f(\xx)\right)\Psi_\sigma(n{\tt{x}}-{\tt{k}})}{\displaystyle\sum_{{\tt{k}}\in \mathcal{A}_n}\Psi_\sigma(n{\tt{x}}-{\tt{k}})}+\frac{\displaystyle\sum_{{\tt{k}}/n\notin R_{ij}^\circ}\left(f\left(\frac{{\tt{k}}}{n}\right)-f(\xx)\right)\Psi_\sigma(n{\tt{x}}-{\tt{k}})}{\displaystyle\sum_{{\tt{k}}\in \mathcal{A}_n}\Psi_\sigma(n{\tt{x}}-{\tt{k}})}\notag\\
&=:T_1^{ij}(\xx)+T_2^{ij}(\xx).
\end{align}
By the convexity of the function $(\cdot)^p$, we have 
$$
A\, \le\ 2^{p-1}\int_{S_{ij}}\left[\left|T_1^{ij}(\xx)\right|^p+\left|T_2^{ij}(\xx)\right|^p\right] d\xx\ =:\ 2^{p-1} \left( A_1 + A_2 \right).
$$
To estimate $A_1$ we proceed similarly as in Theorem \ref{th_nuovo}. In particular, by using Jensen inequality and the first order Taylor formula (\ref{formulataylor}), we get
\begin{equation*}
\begin{split}
A_1\, &\leq\frac{1}{\left[\phis(1)\right]^{2}}\int_{S_{ij}}\sum_{\frac{{\tt{k}}}{n}\in R_{ij}^\circ}\left|f\left(\frac{{\tt{k}}}{n}\right)-f(\xx)\right|^p\Psi_\sigma(n{\tt{x}}-{\tt{k}})\\
&=\frac{1}{\left[\phis(1)\right]^{2}}\int_{S_{ij}}\sum_{\frac{{\tt{k}}}{n}\in R_{ij}^\circ}\left|\sum_{j=1}^2\left(\frac{k_j}{n}-x_j\right)\int_{0}^{1} \frac{\partial f}{\partial x_j}\left({\tt{x}}+t\left(\frac{{\tt{k}}}{n}-{\tt{x}}\right)\right)dt\right|^p\Psi_\sigma(n{\tt{x}}-{\tt{k}})\\
&\leq {M_f\, 2^{p-1} \over \beta_f\, 2\pi} \frac{ M_{p}\left(\phis\right)}{n^{p}\left[\phis(1)\right]^{2}}\int_{S_{ij}}\sum_{j=1}^2\left[\int_0^{2\pi}d\theta\int_{0}^{1}\left|\frac{\partial f}{\partial x_j}\left({\tt{x}}+tC\left(\cos\theta,\sin\theta\right)\right)\right|^p dt\right] d{\tt{x}}\\
&\leq {M_f 2^{p-1}\over \beta_f\, 2\pi} \frac{ M_{p}\left(\phis\right)}{n^{p}\left[\phis(1)\right]^{2}}\int_{I}\sum_{j=1}^2\left[\int_0^{2\pi}d\theta\int_{0}^{1}\left|\frac{\partial f}{\partial x_j}\left({\tt{x}}+tC\left(\cos\theta,\sin\theta\right)\right)\right|^p dt\right] d{\tt{x}}\\
&\leq {M_f \over \beta_f} \frac{K_1^p}{n^{p}}\sum_{j=1}^2\left\|\frac{\partial f}{\partial x_j}\right\|_{p}^{p}.
\end{split}
\end{equation*}
Then, for $A_2$, by the convexity of $|\cdot|^p$ and since each $\Psi_\sigma(n\tt{x}-\tt{k}) \leq 1$, together with Lemma \ref{lemmanuovo} with $\gamma=1$, there exists a constant $c>0$, such that 
\begin{equation*}                        
\begin{split}
A_2\ &\leq\ \frac{1}{\left[\phis(1)\right]^{2p}} \int_{S_{ij}} \sum_{\substack{\frac{\tt{k}}{n} \notin R_{ij}^\circ}} \left| f\left(\frac{\tt{k}}{n}\right) -f(\xx)\right|^p\Psi_\sigma(n\tt{x}-\tt{k})d\tt{x}
       \\
       &\leq \frac{c}{n^{\nu(1-\theta)}\left[\phis(1)\right]^{2p}}\, 2^p\, \|f\|_\infty^p\int_{ S_{ij}}d{\tt{x}}\ \leq\ K_2^p\frac{\|f\|_\infty^p}{n^{\nu(1-\theta)}},
    \end{split}
\end{equation*}
where $0<\nu<\alpha$. Note that Lemma \ref{lemmanuovo} can be applied in the above computations since, being $\xx \in S_{ij}$, $\frac{\tt{k}}{n} \notin R_{ij}^\circ$, it necessarily follows that $\frac{\tt{k}}{n} \in\partial R_{ij}\cup I\setminus R_{ij}$ (where $\partial R_{ij}$ denotes the boundary of $R_{ij}$), and hence $\|n\xx-\kk\| \mau n^{1-\theta}$. 
For what concerns $B$, noting that
$\left|F_n(f,\xx)-f(\xx)\right|^p\miu 2^p\|f\|^p_\infty$, and that the Lebesgue measure of the set $R_{ij}\setminus S_{ij}$ can be bounded by $4\, n^{-\theta}$ we immediately have
$$
B\ \miu 2^{p+2}\|f\|^p_\infty\, n^{-\theta}.
$$
Hence, summing the integrals $\int_{R_{ij}}\left|F_n(f,\xx)-f(\xx)\right|^p\, d\xx$ over all the sets $R_{ij}$ we finally achieve
\begin{equation*}
\begin{split}
    \left\|F_n(f,\cdot)-f\right\|_p^p&\ \leq\ 2^{p-1}\ \sum_{i=0}^{q_1-1}\sum_{j=0}^{q_2-1}\left( {M_f \over \beta_f} \frac{K_1^p}{n^{p}}\left\|f\right\|_{1,p}^p+K_2^p\frac{\|f\|_\infty^p}{n^{\nu(1-\theta)}}+8 \frac{\|f\|_\infty^p}{n^\theta}\right)
    \\
    &=\, 2^{p-1}(q_1\cdot q_2) \left( {M_f \over \beta_f} \frac{K_1^p}{n^{p}}\left\|f\right\|_{1,p}^p+K_2^p\frac{\|f\|_\infty^p}{n^{\nu(1-\theta)}}+8\, \frac{\|f\|_\infty^p}{n^\theta}\right),
    \end{split}
\end{equation*}
from which we get the thesis.
\end{proof}


\section{Error analysis for static greyscale images and implementation of the resizing algorithm} \label{sec2.33}
It is well-known that, a digital image of dimension $N\times M$ (i.e., with $N$ rows and $M$ columns) can be represented as follows
\begin{equation}\label{funz_immagine}
 \widetilde f({\tt{x}})\ :=\ \sum_{i=1}^N\sum_{j=1}^M F_{i,j}\, \chi_{i,j}({\tt{x}}),  \qquad \xx\in I,
\end{equation}
where $\chi_{i,j}({\tt{x}})$ is the characteristic function of the set $R_{ij}$ defined in (\ref{dominio1}) and (\ref{dominio2}) with $l_i=i$, $i=0,...,N$, and $t_j=j$, $j=0,...,M$, and $I=[0,N]\times[0,M]$.
Each element $F_{i,j}$ that for $8$-bit images runs on all integers between 0 and 255, corresponds to the intensity of the pixel at position $(i,j)$ in the greyscale (luminance). Since $\widetilde{f}$ is a piecewice constant function, its partial derivatives are null almost everywhere on the rectangle $[0,N]\times[0,M]$, hence condition (\ref{lower-bound-D}) fails and the results established in the previous sections can not be applied to that case. However we can prove the following.
\begin{theorem}\label{thnuovo3}
Let $\sigma$ be a sigmoidal function satisfying condition $(\Sigma 3)$ with $\alpha>p$, $1 \miu p <+\infty$. Let $\widetilde f$ be defined as in (\ref{funz_immagine}), then we have
\begin{equation*}
\|F_n(\widetilde f,\cdot)-\widetilde f(\cdot)\|_p\ \leq\  255\cdot (2^{p-1}\, MN)^{1/p}\left(\frac{K_2}{n^{\nu(1-\theta)/p}}+\frac{8^{1/p}}{n^{\theta/p}}\right),
\end{equation*}
for $0<\nu<\alpha$ and $0<\theta<1$, where $K_2>0$ arises from Theorem \ref{thnuovo2}.
\end{theorem}
\begin{proof}
Retracing the proof of Theorem \ref{thnuovo2}, in the special case of $\widetilde{f}$, it is easy to see that the corresponding term $A_1$ is identically null, hence the thesis follows taking the estimates of $A_2$ and $B$, and noting that $\|\widetilde{f}\|_\infty \leq 255$, $q_1=N$, $q_2=M$.
\end{proof}
\begin{remark}
The result of Theorem \ref{thnuovo3} can be analogously formulated for any piecewise constant function $f:I\to \R$ (not necessarily in the specific form (\ref{funz_immagine}) which is referred to digital images), as
\be
\|F_n(f,\cdot)-f(\cdot)\|_p\ \leq\ C\, \left(\frac{1}{n^{\nu(1-\theta)/p}}+\frac{1}{n^{\theta/p}}\right),\quad n \in \N,
\ee
for $0<\nu<\alpha$ and $0<\theta<1$, where $C>0$ is a suitable constant.
\end{remark}

Below, we propose an application of the (bivariate) NN operators to digital image reconstruction and resizing. The main idea is to exploit the image representation (\ref{funz_immagine}) to implement a suitable resizing algorithm for magnifying a given digital image by an integer scale factor $r$. It is not restrictive to focus the analysis only on greyscale images since colour images can be naturally handled by processing each RGB channel separately. Later, we also provide specific estimates that can be useful to evaluate the accuracy and the rate of convergence of the proposed algorithm.
\vskip0.2cm

To reach this aim, we have to recall some instances of sigmoidal functions satisfying the above assumptions. 
\vskip0.2cm

A first example is given by the well-known \textit{logistic function} $\sigma_{\ell}(x)=(1+e^{-x})^{-1}$, $x\in\R$, which is very useful in the theory of artificial neural network.  Particularly, in \cite{chen2009approximation} the authors prove that the corresponding density function is given by
\begin{equation*}
\begin{split} 
    \phi_{\sigma_\ell}(x)&=\frac{1}{2}\left(\sigma_\ell(x+1)-\sigma_\ell(x-1)\right)=\frac{e^2-1}{2e^2}\cdot\frac{1}{(1+e^{x-1})(1+e^{-x-1})},\qquad x\in\R.
\end{split}
\end{equation*}
The following useful result holds.
\begin{lemma}[\cite{anastassiou2011multivariate}]\label{logistic}
    Let $0<\theta<1$ and $\gamma>0$. Setting $\Psi_{\sigma_\ell}({\tt{x}})=\prod_{i=1}^d\Phi_{\sigma_\ell}(x_i)$, $d \mau 1$, we have
    \begin{equation*}
        \sum_{{\tt{k}}\in \mathcal{A}_n\atop \|n{\tt{x}}-{\tt{k}}\|>\gamma\, n^{1-\theta}}\Psi_{\sigma_\ell}(n{\tt{x}}-{\tt{k}})\, \leq\, L\, e^{-n^{1-\theta}}, \quad n \in \N,
    \end{equation*}
with $L>0$.
\end{lemma}
Now, we are ready to give the following results in the particular case of NN operators activated by logistic function.
\begin{corollary}\label{corlogistic} 
Under the assumptions of Theorem \ref{thnuovo2}, for $1\miu p <+\infty$, and with $\sigma=\sigma_{\ell}$, we have
\begin{equation*}
\|F_n^{\sigma_\ell}(f,\cdot)-f(\cdot)\|_p\leq \left(2^{p-1}q_1\cdot q_2\right)^{1/p}\, \left({M_f^{1/p} \over \beta_f^{1/p}} \frac{K_1}{n}\left\|f\right\|_{1,p}+8^{1/p}\frac{\|f\|_\infty}{n^{\theta/p}}+K_2\, \|f\|_\infty\, e^{-n^{(1-\theta)}/p}\right),
\end{equation*}
where $F_n^{\sigma_\ell}$ are the NN operators generated by $\Psi_{\sigma_\ell}$, $K_1$, $K_2>0$ are suitable constants and $0<\theta<1$.
\end{corollary}
\begin{proof} 
The proof follows as in Theorem \ref{thnuovo2}, by using Lemma \ref{logistic} in place of Lemma \ref{lemmanuovo}.
\end{proof}

In the special case of digital images $\widetilde{f}$ we are also able to prove a sharper estimate than the general one given in Corollary \ref{corlogistic}; here, exploiting the exponential decay of the function $\phi_{\sigma_\ell}$, we are able to establish the following {\em optimal} order of approximation, by an estimate which turns out to be independent by the technical parameter $\theta$ introduced in the previous results.

\begin{corollary}\label{cor_logistic_imm} Let $\widetilde f$ be defined as in (\ref{funz_immagine}), and $\sigma=\sigma_{\ell}$, then
\begin{equation*}
\|F_n^{\sigma_\ell}(\widetilde f,\cdot)-\widetilde f(\cdot)\|_p\, \leq\,  255 \left(2^{p-1}\, MN\right)^{1/p}C\left[\left(\frac{\ln(n)}{n}\right)^{1/p}+\left(\frac{1}{n}\right)^{1/p}\right] =\, {\cal O}\left({\log(n) \over n}\right),
\end{equation*}
$n \in \N$, $1 \miu p<+\infty$, where $C>0$ is a suitable constant.
\end{corollary}
\begin{proof}
Recalling the notation used in the proof of Theorem \ref{thnuovo2}, we can modify the definition of the sets $S_{ij}$ as follows: 
\[
S_{ij}:=\left[l_i+\frac{\ln(n)}{n},l_{i+1}-\frac{\ln(n)}{n}\right]\times \left[t_i+\frac{\ln(n)}{n},t_{i+1}-\frac{\ln(n)}{n}\right]\subset R_{ij},
\] 
$i=0,\dots,q_1-1$, $j=0,\dots,q_2-1$.
Now, retracing the proof of Theorem \ref{thnuovo2}, we immediately obtain that $A_1=0$. Moreover, concerning the term $A_2$, since $\xx \in S_{ij}$, when $\frac{\tt{k}}{n} \notin R_{ij}^\circ$, it must be that $\frac{\tt{k}}{n} \in \partial R_{ij} \cup (I \setminus R_{ij})$; thus $\left\|n\xx - \kk\right\| \geq \ln(n)$.
Therefore, by the properties of $\phi_{\sigma_\ell}$, we immediately get:
\begin{equation*}
    \begin{split}
\sum_{{\kk \over n} \notin R^{\circ}_{ij}} \Psi_{\sigma_\ell}(n{\tt x}-{\tt k})\ &\leq \sum_{\left\|n{\tt{x}}-{\tt{k}}\right\|> \ln(n)}\Psi_{\sigma_\ell}(n{\tt{x}}-{\tt{k}})\\
&\miu \sum^2_{j=1} \, \left\{ \, \sum_{|nx_j - k_j| > \ln(n)} 
     \phi_{\sigma}(nx_j-k_j) \, \right\}
\\
&\miu \sum^2_{j=1} \, \left\{ \, 2\int_{\ln(n)-1}^{+\infty}\phis(t)dt \, \right\}
\\
&=\ 4 \, \left\{ \, \frac{(-1 + e^{2})\left(\ln\left(n e^{-1} + e^{-1}\right) - \ln\left(1 + n e^{-2}\right) - 1\right)}{ e^{2} (e^{-2} - 1)} \, \right\}\ \miu \ {\widetilde{C} \over n},
\end{split}
\end{equation*}
for $n$ sufficiently large and a suitable $\widetilde{C}>0$. Finally, we also have that $B=\mathcal{O}\left(\ln(n)/n\right)$, as $n\to+\infty$. This completes the proof.
\end{proof}
\begin{remark}
We stress again that, in the specific case considered in Corollary \ref{cor_logistic_imm}, the established estimate is sharper than the more general one of Corollary \ref{corlogistic}.
\end{remark}

As a second interesting example, we can also consider the NN operators activated by the \textit{ramp function} $\sigma_{R}(x)$, defined as
\begin{equation*}
\sigma_{R}(x)=\begin{cases}
0, & \text{ if } x<-\frac{1}{2}, \\
x+\frac{1}{2}, & \text{ if } -\frac{1}{2}\leq x\leq \frac{1}{2}, \\
1, & \text{ if } x>\frac{1}{2}.
\end{cases}
\end{equation*} 
The corresponding density function for the ramp activation, $\phi_{\sigma_{R}}$, can be expressed in terms of ReLU (\textit{Rectified Linear Unit}) activation function, defined by $\psi_{\text{ReLU}}(x)=(x)_{+}:= \max \{x, 0\}$, $x\in\R$. The expression for $\phi_{\sigma_{R}}$ is given by
\begin{equation*}
\phi_{\sigma_R}(x)=\frac{1}{2}\left[\psi_{\text{ReLU}}\left(x+\frac{3}{2}\right)-\psi_{\text{ReLU}}\left(x+\frac{1}{2}\right)-\psi_{\text{ReLU}}\left(x-\frac{1}{2}\right)+\psi_{\text{ReLU}}\left(x-\frac{3}{2}\right)\right],
\end{equation*}
$x\in\R$. For further details, we refer to \cite{daubechies2022nonlinear, piconi2024, coroianu2024approximation, opschoor2022exponential, li2024two, CJat23,chen2022construction,plonka2023spline}.\\

For the special case of the ramp function we can state the following corollary, that for brevity we directly write in the case of functions $\widetilde{f}$.

\begin{corollary} \label{corollary-ramp-function}
Let $\widetilde f$ be defined as in (\ref{funz_immagine}), and $\sigma=\sigma_{R}$. Then
\begin{equation*}
\|F_n^{\sigma_R}(\widetilde f,\cdot)-\widetilde f(\cdot)\|_p\, \leq\, 255 \left(2^{p-1}\, MN\right)^{1/p}\, C\,\left(\frac{\ln(n)}{n}\right)^{1/p}, 
\end{equation*}
$1 \miu p<+\infty$, for sufficiently large $n \in \N$, $C>0$.
\end{corollary}
\begin{proof}
The proof is the same of Corollary \ref{cor_logistic_imm}; here we only have to observe that, since $\Psi_{\sigma_R}$ has compact support, for $n$ sufficiently large also $A_2=0$. This completes the proof.
\end{proof}

In the next sections, we provide numerical experiments concerning the implementation of the procedure for image reconstruction and resizing.
The reconstruction method has been implemented in \textsc{Matlab\textsuperscript{\textcopyright}}; for completeness, below we provide the pseudocode of the algorithm based on NN operator, that we call in what follows (for brevity) \textit{NN algorithm}; see, Algorithm \ref{alg}.

\RestyleAlgo{ruled}
\begin{algorithm}[hbt!]
\caption{{\small Pseudocode of the NN algorithm}}\label{alg}
\textbf{Goal:} Reconstructing and enhancing (resizing) the resolution of the original image $I$ using NN operators activated by the logistic or ramp functions.\\
\textbf{Input data:} Original image $I$ ($N\times M$ pixel resolution); parameter $n\in\mathbb{N}^+$ of the operator $F_n$; integer scaling factor $r>1$.\\
- Choice and definition of the density functions $\phi_{\sigma}$ and $\Psi_{\sigma}$\;
- Computation of the size of the reconstructed image using the scaling factor: $(N\cdot r)\times (M\cdot r)$\;
- Modelling of $I$ as a function $\widetilde{f}_I$ of the form (\ref{funz_immagine})\;
- Definition of a sampling grid, $G_{r}$, containing $(N\cdot r)\times (M\cdot r)$ uniformly spaced nodes over the square domain $[0,N]\times[0,M]$\;
- Computation of the matrix of the sample values $f({\tt{k}}/n)$, ${\tt{k}}\in\mathcal{A}_n$\;
- Definition of the vectors containing the arguments of $\Psi_{\sigma}$\;
 \For{$i=1,\dots,N$ and $j=1,\dots,M$}{
  sum over ${\tt{k}}\in\mathcal{A}_n$ of all non-zero or not negligible terms of the form $\Psi_\sigma(n{\tt{x}}_{ij}-{\tt{k}})f\left({\tt{k}}/n\right)$ for ${\tt{x}}_{ij}$ belonging to the reconstruction grid}\;
\KwResult{The reconstructed image with increased resolution $(N\cdot r)\times (M\cdot r)$ consisting in the evaluation of $F_n \widetilde{f}_I$ at the points of the grid $G_{r}$}.
\end{algorithm}
The pseudocode in Algorithm \ref{alg} relies on the idea of reconstructing a resized image (by a scaling factor $r$) by evaluating the operator $F_n \widetilde{f}_I$ over the points of the sampling grid $G_r$. The number of nodes in $G_r$ depends by the scaling factor $r>1$, and it is equal to $(N \cdot r) \times (M \cdot r)$ which corresponds to the final size of the rescaled image. A fundamental step in the algorithm is the definition of the points of $G_r$ that must be chosen avoiding the borders of the sets $R_{ij}$ where $\widetilde{f}_I$ is discontinuous.
%

\section{Rate of convergence in $L^p$-spaces via a constructive density method by the averaged modulus of smoothness} \label{secLp}

In what follows, we extend the estimates provided in the previous section to the whole $L^p$-spaces. In order to do this, we recall the following useful tool, that is the so-called (multivariate) averaged modulus of smoothness (or $\tau$-modulus), \cite{popovKH1983averaged}.
\vskip0.2cm

For every $f \in L^p(I)$ we first recall the usual (first-order) modulus of smoothness
$$
\omega(f,\delta)_p\ :=\ \sup_{0<|{\tt h}|\miu\delta} \left(   \int_I \left|   f(\xx+{\tt h}) - f(\xx) \right|^p\, d\xx \right)^{1/p}, \quad 1 \miu p<+\infty.
$$

Now, denoting by $M(I)$ the set of all bounded and measurable functions $f:I\to \R$, we can define the (multivariate) local modulus of smoothness as
$$
\omega(f,\xx;\, \delta)\ :=\ \sup_{0<|{\tt h}|\miu\delta} \left\{    \left|   f(\uu+{\tt h}) - f(\uu) \right|:\ \uu+{\tt h},\, \uu \in \Omega_{\delta/2}(\xx) \right\},  \quad \xx \in I,
$$
where
$$
\Omega_{\delta/2}(\xx)\ :=\ \left\{ \yy \in \R^2:\ \|\yy-\xx\|_m \miu \delta/2 \right\}.
$$

The (multivariate) averaged modulus of smoothness (or first order $\tau$-modulus) in the $L^p$-norm can be now defined as
\be
\tau(f, \delta)_p\ :=\ \| \omega(f,\cdot;\, \delta)\|_p, \quad 1 \miu p <+\infty.
\ee
Note that, the $\tau$-modulus can be defined also for $p=+\infty$ when $f \in C(I)$, and it turns out that
$$
\tau(f,\delta)_\infty\ =\ \omega(f,\delta), \quad \delta>0,
$$
where $\omega(f,\delta)$ is the classical modulus of continuity of $f$ (see, e.g., \cite{devore1993constructive}).\\

Below we recall some important and useful properties of the $\tau$-modulus (for more details, see \cite{popovKH1983averaged}).
\begin{lemma} \label{lemma-tau}
Let $f$, $g \in M(I)$, and $1 \miu p<+\infty$. Thus
\begin{itemize}

\item[(i)] $\tau(f,\delta)_p \miu \tau(f, \delta')_p$, for $0<\delta \miu \delta'$;

\item[(ii)] $\tau(f+g, \delta)_p\miu \tau(f,\delta)_p+\tau(g,\delta)_p$;

\item[(iii)] $\tau(f,\lambda\, \delta)_p \miu (2\lfloor \lambda\rfloor +2)^{3} \tau(f,\delta)_p$, $\lambda>0$;

\item[(iv)] $\displaystyle \tau(f,\delta)_p\, \miu\, 2\, \delta\, \sum_{i=1}^2 \left\| {\partial f \over \partial x_i} \right\|_p$, with $f\in W^{1}_{p}(I)$;

\item[(v)] $\omega(f,\delta)_p\ \miu\ \tau(f,\delta)_p$.

\end{itemize}
\end{lemma}

In order to achieve the desired result we can prove the following proposition.
\begin{proposition} \label{prop1}
For every $f \in M(I)$ there exists piecewice constant functions $P_f$, $Q_f \in M(I)$ such that the following estimates holds
\begin{itemize}

\item[(A)] $Q_f \miu f \miu P_f$;

\item[(B)] $\|f - P_f\|_p \miu C_1\, \tau(f, n^{-1/4})_p$, $\|f - Q_f\|_p \miu C_2\, \tau(f, n^{-1/4})_p$, $n \in \N$ sufficiently large;

\item[(C)] for any fixed $\sigma$ satisfying $(\Sigma 3)$ for $\alpha>1$, we have: $\|F_n(P_f,\cdot) - P_f(\cdot)\|_1 \miu C_3\, \tau(f, n^{-1/4})_1$, $\|F_n(Q_f, \cdot) - Q_f(\cdot)\|_1 \miu C_4\, \tau(f, n^{-1/4})_1$, $n \in \N$ sufficiently large;

\end{itemize}
for suitable constants $C_i>0$, $i=1,...,4$.
\end{proposition}
\begin{proof}
The proof of such proposition is constructive. Let $n \in \N$ sufficiently large, and we consider $n_1=n_2:=\lfloor n^{1/4} + 1\rfloor$ such that $h_1:=(b_1-a_1)/n_1 \miu (b_1-a_1)n^{-1/4}$, and $h_2:=(b_2-a_2)/n_2 \miu (b_2-a_2)n^{-1/4}$. Using $h_1$ and $h_2$ we can consider a decomposition of $I$ with rectangles $R_{ij}$ as in (\ref{dominio1}) and (\ref{dominio2}), defined with nodes
$$
l_i:= a_1+i\, h_1, \quad t_j:=a_2+j\, h_2, \quad i = 0, ..., n_1,\ j=0, ..., n_2.
$$
Now, denoting by $R^*_{ij}$ the closure of $R_{ij}$, we can define the following real values
$$
P_{ij}:=\sup\left\{ f(\xx):\ \xx \in R^*_{ij} \right\}, \quad Q_{ij}:=\inf\left\{ f(\xx):\ \xx \in R^*_{ij} \right\},
$$
for $i = 0, ..., n_1$, $j=0, ..., n_2$. Using these values we can define the following two piecewise constant functions
\be \label{auxiliary-functions}
P_f(\xx):=P_{ij},\ \text{if}\ \xx \in R_{ij}, \quad \quad Q_f(\xx):=Q_{ij},\ \text{if}\ \xx \in R_{ij}.
\ee
It is clear by the above construction that $P_f$ and $Q_f$ belong to $M(I)$ and satisfy $(A)$.\\
 We now prove $(B)$.
Observing that the inclusion
$$
R^*_{ij} \subset \Omega_{|{\tt h}|}(\xx), \quad \xx \in R_{ij}, \quad {\tt h}:=(h_1,h_2),
$$
holds, for $i = 0, ..., n_1$, $j=0, ..., n_2$, we can get the following inequality
$$
\left| P_f(\xx)-f(\xx)  \right|\, \miu\, \omega(f, \xx;\, 2|{\tt h}|), \quad \xx \in I.
$$
Indeed, for $\xx \in R_{ij} \subset I$, we have
$$
\left| P_f(\xx)-f(\xx)  \right|\, =\, \left( \sup_{{\tt y} \in R^*_{ij}}  f({\tt y}) \right) - f(\xx)\ \miu\, \sup_{{\tt y},\, {\tt z} \in R^*_{ij}}  [f({\tt y})-f({\tt z})]
$$
\be \label{stima-tau-doppio}
\miu\, \sup_{{\tt y},\, {\tt z} \in \Omega_{|{\tt h}|}(\xx)}  |f({\tt y})-f({\tt z})|\ =\ \omega(f, \xx;\, 2|{\tt h}|).
\ee
Passing to the $p$-power of the previous inequality and integrating over $I$, we get
$$
\| P_f-f\|_p\, \miu\ \tau(f, 2|{\tt h}|)_p.
$$
Setting $T:=\max\{ b_1-a_1,\ b_2-a_2\}$, and using Lemma \ref{lemma-tau} (i) and (iii), we finally obtain
$$
\| P_f-f\|_p\, \miu\ \tau(f, 4\, T n^{-1/4})_p\ \miu\ C_1\, \tau(f, n^{-1/4})_p,
$$
for a suitable $C_1>0$. Since the second part of $(B)$ can be proved analogously, we directly turn to the proof of $(C)$. We define the following auxiliary functions $\Theta_{ij}:I \to \R$
$$
\Theta_{ij}(\xx)\, =\, \Theta_{ij}(x_1,x_2)\ :=\, \begin{cases}
0, \quad \text{if } x_1 <l_i\, \text{ or } \, x_2<t_j,\\
1,  \quad \text{if } x_1 \mau l_i\, \text{ and }\, x_2\mau t_j,
\end{cases}
$$
$i,\, j = 0, ..., n$. Using $\Theta_{ij}$ we can write the following functional representation for the step function $P_f$
$$
P_f(\xx)\, =\, \sum_{i=0}^{n_1-1}\sum_{j=0}^{n_2-1} \left\{   P_{ij}-P_{i\ j-1} - P_{i-1\ j}+P_{i-1\ j-1}\right\}\, \Theta_{ij}(\xx), \quad \xx \in I,
$$ 
where in the above formula we put $P_{i\ -1}=P_{-1\ j}=0$, $i = 0, ..., n_1$, $j=0, ..., n_2$. Indeed, for any $\xx \in R_{\mu \nu}$ one has
$$
\sum_{i=0}^{n_1-1}\sum_{j=0}^{n_2-1} \left\{   P_{ij}-P_{i\ j-1} - P_{i-1\ j}+P_{i-1\ j-1}\right\}\, \Theta_{ij}(\xx)\ 
$$
$$
=\, \sum_{i=0}^{\mu}\sum_{j=0}^{\nu} \left\{   P_{ij}-P_{i\ j-1} - P_{i-1\ j}+P_{i-1\ j-1}\right\}\, =\, \sum_{i=0}^{\mu}\left[\sum_{j=0}^{\nu} \left\{   P_{ij}-P_{i\ j-1} \right\}\right]\,
$$
$$
  -\, \sum_{i=0}^{\mu}\left[\sum_{j=0}^{\nu} \left\{   P_{i-1\ j}-P_{i-1\ j-1}\right\}\right]\, 
=\, \sum_{i=0}^{\mu} P_{i\, \nu} - \sum_{i=0}^{\mu} P_{i-1\, \nu}\ =\ P_{\mu\, \nu}\ =\ P_f(\xx).
$$
Now, recalling the above expression for $P_f$, the linearity of the operators $F_n$, and the Minkowsky inequality, one can write what follows
$$
\|F_n(P_f, \cdot)-P_f(\cdot)\|_1\ \miu\ \sum_{i=0}^{n_1-1}\sum_{j=0}^{n_2-1} \left|  P_{ij}-P_{i\ j-1} - P_{i-1\ j}+P_{i-1\ j-1}\right|\, \cdot \|F_n(\Theta_{ij}, \cdot)-\Theta_{ij}(\cdot)\|_1.
$$
Proceeding similarly as in (\ref{stima-tau-doppio}), we can observe that
$$
\left|  P_{ij}-P_{i\ j-1} - P_{i-1\ j}+P_{i-1\ j-1}\right|\, \miu\, 2\, \omega(f, {\tt n}_{ij}; 2|{\tt h}|), 
$$
where ${\tt n}_{ij}:=(l_i,t_j)$, $i= 0, ..., n_1-1$, $j=0,...,n_2-1$, obtaining
\be \label{ggtthhqq}
\|F_n(P_f, \cdot)-P_f(\cdot)\|_1\ \miu\ \sum_{i=0}^{n_1-1}\sum_{j=0}^{n_2-1}\, \omega(f, {\tt n}_{ij}; 2|{\tt h}|)\, \cdot \|F_n(\Theta_{ij}, \cdot)-\Theta_{ij}(\cdot)\|_1.
\ee
Now, we have to estimate the terms $\|F_n(\Theta_{ij}, \cdot)-\Theta_{ij}(\cdot)\|_1$. Proceeding as in the proof of Theorem \ref{thnuovo2} and Theorem \ref{thnuovo3} we have
$$
\|F_n(\Theta_{ij}, \cdot)-\Theta_{ij}(\cdot)\|_1
$$
$$
\leq\ \left\{  \int_{l_i}^{b_1}\int_{t_j}^{b_2} + \int_{a_1}^{l_{i}-{1/\sqrt{n}}} \int_{a_2}^{b_2} + \int_{l_{i}-1/\sqrt{n}}^{l_{i}} \int_{a_2}^{b_2} +     \int_{a_1}^{b_1} \int_{a_2}^{t_{j}-1/\sqrt{n}} +   \int_{a_1}^{b_1} \int_{t_{j}-1/\sqrt{n}}^{t_j} \right\} \left|  F_n(\Theta_{ij}, \xx)-\Theta_{ij}(\xx) \right|\, d\xx
$$
$$
=:\ S_1\ +\ S_2\ +\ S_3\ +\ S_4\ + S_5\footnote{We here represented only the more general situation of indexes $i\, j$; the particular cases arising when $i$ is equal to $0, 1, n_1-1$ or $j$ is equal to $0,1,n_2-1$ can be treated splitting the integral with less term than the above general case.}.
$$
By the definition of $\Theta_{ij}$, we have
\begin{equation*}
\begin{split}
S_2\ &=\ \int_{a_1}^{l_{i}-1/\sqrt{n}} \int_{a_2}^{b_2} \left|  F_n(\Theta_{ij}, \xx) \right|\, d\xx\ \miu\ {1 \over (\phis(1))^2} \int_{a_1}^{l_{i}-1/\sqrt{n}} \int_{a_2}^{b_2} \sum_{\kk \in {\cal A}_n} \left| \Theta_{ij}\left( {k_1 \over n}, {k_2 \over n}\right) \right| \Psi_{\sigma}(n \xx-\kk)\, d\xx\\
&=\ {1 \over (\phis(1))^2} \int_{a_1}^{l_{i}-1/\sqrt{n}}  \sum_{k_1/n \mau l_i} \phis(nx_1-k_1)\, dx_1 \int_{a_2}^{b_2} \sum_{k_2/n \mau t_j} \phis(nx_2-k_2)\, dx_2\\
&\miu {(b_2-a_2) \over (\phis(1))^2} \int_{a_1}^{l_{i}-1/\sqrt{n}}  \sum_{k_1/n \mau l_i} \phis(nx_1-k_1)\, dx_1\ \miu\ {T^2 \over (\phis(1))^2} \sup_{x \in \R} \sum_{|nx-k_1|>\sqrt{n}} \phis(nx-k_1) \miu K_1 n^{-1/2},
\end{split}
\end{equation*}
for $n$ sufficiently large, by Lemma \ref{lemmanuovo} with $\nu=1$, $\theta=1/2$, and for a suitable $K_1>0$. The same estimate holds also for $S_4$. By the same reasoning, we can immediately obtain
$$
S_k\ \miu\ {T \over \sqrt{n}}, \quad \mbox{for} \quad k=3, 5.
$$
 It remains to estimate $S1$. In this case, writing (\ref{spezzato}) with the function $\Theta_{ij}$, we obtain
\begin{equation*}
\begin{split}
&\int_{l_i}^{b_1}\int_{t_j}^{b_2} |F_n(\Theta_{ij}, \xx)-\Theta_{ij}(\xx)|\, d\xx\\
&\miu\ {1 \over (\phis(1))^2} \left\{ \int_{l_i}^{b_1}\int_{t_j}^{b_2} \sum_{\lceil na_1  \rceil \miu k_1<n l_i}\ \sum_{k_2=\lceil na_2  \rceil}^{\lfloor nb_2 \rfloor}  \Psi_\sigma(n \xx -\kk)\ d\xx \right.\\
&\left. + \ \int_{l_i}^{b_1}\int_{t_j}^{b_2} \sum_{k_1=\lceil na_1  \rceil}^{\lfloor nb_1 \rfloor} \sum_{\lceil na_1  \rceil \miu k_2<n t_j}\ \Psi_\sigma(n \xx -\kk) d\xx \right\}\ =:\ {1 \over (\phis(1))^2} \left\{ S_{1,1}+S_{1,2}\right\}.
\end{split}
\end{equation*}
For $S_{1,1}$ we have
\begin{equation*}
\begin{split}
S_{1,1} &\miu (b_2-t_j)\, \left\{ \int_{l_i}^{l_i+1/\sqrt{n}} \sum_{\lceil na_1  \rceil \miu k_1<l_i n} \phis(nx_1-k_1)\ dx_1 +\ \int_{l_i+1/\sqrt{n}}^{b_1} \sum_{\lceil na_1  \rceil \miu k_1<n l_i} \phis(nx_1-k_1) dx_1 \right\}\\
&\miu T\, \left\{ {1 \over \sqrt{n}} +\ \int_{l_i+1/\sqrt{n}}^{b_1} \sum_{\lceil na_1  \rceil \miu k_1<n l_i} \phis(nx_1-k_1) dx_1 \right\} \miu {K_2 \over \sqrt{n}},
\end{split}
\end{equation*}
for a suitable $K_2>0$, where the last part of the estimate follows by the same procedure adopted for estimating $S_2$, and using Lemma \ref{lemmanuovo} with $\nu=1$, $\theta=1/2$. Finally, the term $S_{1,2}$ can be estimated analogously. Rearranging all the above terms, we finally have with a suitable $C>0$
$$
\|F_n(\Theta_{ij}, \cdot)-\Theta_{ij}(\cdot)\|_1\ \miu C\, {1 \over \sqrt{n}}.
$$
Using the above estimate in (\ref{ggtthhqq})
$$
\|F_n(P_f, \cdot)-P_f(\cdot)\|_1\ \miu\  C\, {1 \over \sqrt{n}}  \sum_{i=0}^{n_1-1}\sum_{j=0}^{n_2-1}\,  \omega(f, {\tt n}_{ij}; 2|{\tt h}|) 
$$
$$
=\ C\, {1 \over h_1\, h_2\, \sqrt{n}}  \sum_{i=0}^{n_1-1}\sum_{j=0}^{n_2-1}\, \int_{R_{ij}} \omega(f, {\tt n}_{ij}; 2|{\tt h}|)\, d{\tt x}\ \miu\ \widetilde{C}\, \sum_{i=0}^{n_1-1}\sum_{j=0}^{n_2-1}\, \int_{R_{ij}} \omega(f, {\tt n}_{ij}; 2|{\tt h}|)\, d{\tt x},
$$
for every sufficiently large $n$. Now, observing that
$$
\Omega_{|{\tt h}|}({\tt n}_{ij})\ \subset \Omega_{2|{\tt h}|}(\xx), \quad \xx \in R_{ij},
$$
we immediately obtain
$$
\omega(f, {\tt n}_{ij}; 2|{\tt h}|)\ \miu\ \omega(f, \xx; 4|{\tt h}|), \quad \xx \in R_{ij},
$$
and thus
$$
\|F_n(P_f, \cdot)-P_f(\cdot)\|_1 \miu \widetilde{C}\, \sum_{i=0}^{n_1-1}\sum_{j=0}^{n_2-1}\, \int_{R_{ij}} \omega(f, \xx; 4|{\tt h}|)\, d{\tt x} = \widetilde{C}\, \int_{I} \omega(f, \xx; 4|{\tt h}|)\, d{\tt x} = \widetilde{C} \tau(f, 4|{\tt h}|)_1.
$$
Hence, by using Lemma \ref{lemma-tau} (i) and (iii), and the definition of $h_1$ and $h_2$, immediately get
$$
\|F_n(P_f, \cdot)-P_f(\cdot)\|_1 \miu C_3\, \tau(f, n^{-1/4})_1
$$
for a suitable $C_3>0$. This completes the proof of the first part of $(C)$; clearly the second part involving $Q_f$ follows by the same reasoning.
\end{proof}

Now, in order to establish the main result of this section, we have to recall the following interpolation theorem of Riesz-Thorin type involving the $\tau$-modulus and positive linear operators. This result can be found, e.g., in Theorem 5.1 of \cite{quak1989multivariate}.
\begin{theorem}[Riesz-Thorin type theorem] \label{Riesz-Thorin}
Let $L:M(I) \to M(I)$ be a given positive linear operator, such that for $\delta>0$ there holds
$$
\|Lf-f\|_{\widetilde{p}}\ \miu\ M_1\, \tau(f,\delta)_{\widetilde{p}}, \quad f \in M(I),
$$
$$
\|Lf-f\|_{p*}\ \miu\ M_2\, \tau(f,\delta)_{p*}, \quad f \in M(I),
$$
with $1 \miu \widetilde{p},\, p* \miu +\infty$, and for positive absolute constant $M_1$ and $M_2$. Then, for every $f \in M(I)$ the following estimate holds true
$$
\|Lf -f\| \miu M_3\, \tau(f,\delta)_{p},
$$
for every $1 \miu p \miu +\infty$, such that
$$
{1 \over p}\ =\ {1-\theta \over \widetilde{p}}\ +\ {\theta \over p*}, \quad 0<\theta<1.
$$
\end{theorem}
Now we are able to prove the following.
\begin{theorem} \label{Lp-estimates}
Let $\sigma$ be a sigmoidal functions satisfying $(\Sigma 3)$ for $\alpha>1$. Then, for every $f \in M(I)$ the following estimate there holds
$$
\|F_n(f, \cdot)-f(\cdot)\|_p\ \miu\ C\, \tau(f,n^{-1/4})_p, \quad 1 \miu p \miu +\infty,
$$
$n \in \N$ sufficiently large, where $C>0$ is a suitable constant.
\end{theorem}
\begin{proof}
We first prove the thesis for $p=1$. For a given $f \in M(I)$, using Proposition \ref{prop1} (B) and (C) we have
\begin{equation*}
\begin{split}
\|f(\cdot) - F_n(f, \cdot)\|_1\, &\miu\, \|f(\cdot)-P_f(\cdot)\|_1\, +\, \|P_f(\cdot)-F_n(P_f, \cdot)\|_1 +\, \|F_n(P_f, \cdot)-F_n(f,\cdot)\|_1\\
&\miu\ C_1\, \tau(f, n^{-1/4})_1\, +\, C_3\, \tau(f, n^{-1/4})_1\, +\, \|F_n(P_f, \cdot)-F_n(f,\cdot)\|_1.
\end{split}
\end{equation*}

Moreover, for the last term, using Proposition \ref{prop1} (A), and subsequently Proposition \ref{prop1} (B) and (C) again, we can finally write
\begin{equation*}
\begin{split}
&\|F_n(P_f, \cdot)-F_n(f,\cdot)\|_1\ \miu\  \|F_n(P_f, \cdot)-F_n(Q_f,\cdot)\|_1\\
&\miu\ \|F_n(P_f, \cdot)-P_f(\cdot)\|_1 + \|P_f(\cdot)-f(\cdot)\|_1 + \|f(\cdot)-Q_f(\cdot)\|_1 + \|Q_f(\cdot)-F_n(Q_f,\cdot)\|_1\\
&\miu\ \left( C_1+C_2+C_3+C_4 \right) \tau(f, n^{-1/4})_1,
\end{split}
\end{equation*}
hence the proof for the case $p=1$ immediately follows with a suitable constant $C>0$.
\vskip0.2cm

Concerning the case $p=+\infty$, noting that $\tau(f, \delta)_\infty=\omega(f, \delta)$, $\delta>0$, the proof immediately follows by Theorem 11, 12 and 13 of \cite{coroianu2022quantitative}.
\vskip0.2cm

Finally, all the cases of $1<p<+\infty$ follow by the application of Theorem \ref{Riesz-Thorin}.
\end{proof}

\begin{remark} \label{remark-ultima}
Note that the result established in Theorem \ref{Lp-estimates} implies the $L^p$-convergence of the NN operators $F_n$ for any Riemann integrable function $f \in M(I)$ (see Proposition 4 of \cite{bardaro2014generalized}). 

Furthermore, the estimate established in Theorem \ref{Lp-estimates} is more general than the one achieved in Theorem \ref{th_nuovo}, Theorem \ref{thnuovo2}, and Theorem \ref{thnuovo3}, however it is less sharp than the previous ones. This consideration always holds when comparing Theorem \ref{Lp-estimates} with Theorem \ref{th_nuovo}, while with the other two theorems is valid, e.g., if we set $\nu=p$, $\theta=1/2$, and $1 \miu p< 2$.
\end{remark}

For the sake of completeness, we can also state the following corollary as a consequence of Lemma \ref{lemma-tau} (iv).

\begin{corollary} \label{ultimo-cor}
Let $\sigma$ be a sigmoidal functions satisfying $(\Sigma 3)$ for $\alpha>1$. Then, for every Riemann integrable and bounded $f \in W^1_p(I)$ the following estimate there holds
$$
\|F_n(f, \cdot)-f(\cdot)\|_p\ \miu\ C\, n^{-1/4} \sum_{i=1}^2\left\| {\partial f \over \partial x_i}  \right\|_p, \quad 1 \miu p \miu +\infty,
$$
$n \in \N$ sufficiently large, where $C>0$ is a suitable constant.
\end{corollary} 

Concerning Corollary \ref{ultimo-cor}, we can give similar considerations to those given in Remark \ref{remark-ultima} in relation to the achieved order of approximation and its comparison with the result of Section \ref{sec2.2} and Section \ref{sec2.33}. Finally, estimates for the NN operators with respect to the averaged modulus of smoothness, for functions of one-variable, can be found in \cite{CJat23}.



\section{Convergence rate in terms of SSIM and cSSIM indexes}\label{3.1}

The \textit{structural similarity index} (SSIM) \cite{UIQIindex, SSIM1, SSIM2} is an indicator used to measure the similarity of two images (e.g., in the presence of noise or distortion, or after the application of a given algorithm). Rather than focusing on individual pixel differences, SSIM compares the structural patterns of pixel intensities between neighboring pixels, after normalizing for luminance and contrast variations. In the last years, its mathematical properties have been extensively studied, and different results have been derived (see, e.g., \cite{brunet2012geodesics, brunet2011mathematical}).\\
Given two images $F,G\in \R^{N\times M}_{\geq 0}$ (where $\R^{N\times M}_{\geq 0}$ denotes the set of all $N\times M$ matrices with non-negative real entries), with $N,M\in\mathbb{N}$ representing the dimensions of the images, the SSIM is defined by
\begin{equation*}\label{def_SSIM}
\text{SSIM}(F,G) := l(F,G)\cdot c(F,G)\cdot s(F,G)= \frac{2\mu_F\mu_G + C_1}{\mu_F^2 + \mu_G^2 + C_1}\cdot\frac{2\sigma_{FG} + C_2}{\sigma_F^2 + \sigma_G^2 + C_2},
\end{equation*}
where $l(F,G) = \frac{2\mu_F\mu_G + C_1}{\mu_F^2 + \mu_G^2 + C_1}$ measures the difference in luminance between two images, and $c(F,G)=\frac{2\sigma_F\sigma_G + C_1}{\sigma_F^2 + \sigma_G^2 + C_1}$ measures the difference in contrast. Here $\mu_F$, $\mu_G$ are the sample mean of $F$ and $G$, while $\sigma_F^2$ , $\sigma_G^2$, and $\sigma_{FG}$ are the sample variances and covariance, respectively. The structure component $s(F,G)$ evaluates the correlation of local structures in images, considering the similarity of patterns, textures and details. Hence, an SSIM value of $1$ indicates perfect similarity, while a value of $0$ indicates no similarity. The constants $C_1, C_2>0$ are small stabilizing factors that prevent division by zero, and they are usually set as  $C_1 = (0.01 \times L)$, $C_2 =(0.03 \times L)$, where $L=255$ denotes the dynamic range of pixel values for 8-bit images.\\

Recently, a continuous extension of SSIM, known as the \textit{continuous SSIM} (cSSIM), has been introduced in \cite{marchetti2021}. Such formulation extends the discrete SSIM to a continuous framework, which allows one to analyze sequences of images of increasingly fine resolutions.\\
Let $\Omega\subset \R^d$ be a bounded domain, and $\nu$ be a probability measure on $\Omega$. Denote by $L^2_{+}(\Omega)$ the set of functions $f\in L^2(\Omega)$ that are non-negative $\nu$-almost everywhere. For $f,g\in L^2_+(\Omega)$, the cSSIM is defined as follows
\begin{equation*}\label{def_cSSIM}
\text{cSSIM}(f,g) := \frac{2\mu_f\mu_g + c_1}{\mu_f^2 + \mu_g^2 + c_1}\cdot\frac{2\sigma_{fg} + c_2}{\sigma_f^2 + \sigma_g^2 + c_2},
\end{equation*}
with the constants $c_1, c_2 > 0$ have the same stabilizing role as in the discrete case. Here
\begin{equation*}
\begin{split}
\mu_f &:=\int_{\Omega}fd\nu=\left\|f\right\|_1,\\
\sigma_{fg} &:=\int_{\Omega}(f-\mu_f)(f-\mu_g)d\nu,\\
\sigma_f^2 &:= \sigma_{ff}.
\end{split}
\end{equation*}
\begin{remark}
The SSIM can be seen as a special case of cSSIM when the measure $\nu$ is taken as the discrete counting measure. Conversely, the discrete SSIM is essentially a discretized version of cSSIM. As continuous functions $f,g$ are discretized into images with progressively finer resolutions, the SSIM of these discretizations converges to the cSSIM of the original functions. This relationship is formalized in Proposition 2.2 of \cite{marchetti_santin2022}
\end{remark}

The following result shows that the approximation error between two functions $f$ and $g$ by the means of the cSSIM can be controlled by the $L^2$-norm distance of the two functions. 
\begin{theorem}[Theorem 3.3 of \cite{marchetti_santin2022}]\label{th_marchettisantin}
Let $f,g$ be two positive real valued functions in $L^2(I)$ and $c_1,c_2>0$. Then, the following bound holds
{\rm
\begin{equation*}
|1-\text{cSSIM}(f,g)|\leq c_f \left\|f-g\right\|_2^2,
\end{equation*}}
where $c_f:=\frac{4}{\sigma_f^2+c_2}+\frac{1}{\mu_f^2+c_1}$.
\end{theorem}
The metric $|1-\text{cSSIM}(f,g)|$ is called the (continuous) {\em dissimilarity index} between $f$ and $g$. Since we already discussed the importance of SSIM (and, of course, its continuous version) in image analysis, and after the applications discussed in Section \ref{sec2.33}, it is natural to ask if for the considered family of linear operators, convergence results with respect to the continuous dissimilarity index can be proven, when functions like those in (\ref{funz_immagine}) are considered.
\begin{theorem}\label{corollary_nuovo} 
Let $0<\theta<1$ be fixed, and $\sigma$ be a sigmoidal function satisfying $(\Sigma 3)$ with $\alpha>\theta/(1-\theta) \geq 2$. For any function $\widetilde{f}$ of the form (\ref{funz_immagine}), there holds 
\begin{equation*}
|1-\text{cSSIM}(\widetilde{f}(\cdot),\, F_n(\widetilde{f}, \cdot))|\ \leq\ c_f\, 255^2 M\, N\, \frac{K}{n^{\theta}}\,,  
\end{equation*}
for $n \in \N$ sufficiently large, a suitable $K>0$, and where $c_f$ arises from Theorem \ref{th_marchettisantin}.
\end{theorem}
\begin{proof}
Let $\widetilde{f}$ be as in the statement. By Theorem \ref{thnuovo3} with $\nu =\theta/(1-\theta)$ we immediately obtain
$$
\| F_n(\widetilde{f}, \cdot) - f(\cdot) \|_2^2\, \miu\ 255^2 M\, N\, \frac{K}{n^{\theta}},
$$
where $K>0$ is a suitable constant, and $n$ is sufficiently large. Thus, applying Theorem \ref{th_marchettisantin}, we immediately get the thesis.
\end{proof}
Again, Theorem \ref{corollary_nuovo} is valid for general sigmoidal functions $\sigma$; however, as observed in Section \ref{sec2.33}, in the specific cases of the logistic and ramp functions, we can obtain a sharper convergence rate also with respect to the dissimilarity index.
\begin{corollary}\label{CorollaryFinal}
For any function $\widetilde{f}$ of the form (\ref{funz_immagine}) and $\sigma=\sigma_R$ or $\sigma=\sigma_\ell$, denoting by $F_n^\sigma$ the operator generated by $\sigma$, there holds 
\begin{equation*}
|1-\text{cSSIM}(\widetilde{f}(\cdot),\, F_n^{\sigma}(\widetilde{f}, \cdot))|\ \leq\ c_f\, 255^2 M\, N\, K\, \frac{\ln(n)}{n}\,,  
\end{equation*}
for $n\in\N$ sufficiently large, where $c_f$ is the constant arising from Theorem \ref{th_marchettisantin} and $K>0$ is a suitable constant.
\end{corollary}

It should be noted again that the order established in Corollary \ref{CorollaryFinal} (that is an immediate consequence of Corollary \ref{cor_logistic_imm} and Corollary \ref{corollary-ramp-function}) does not depend on the technical parameter $\theta$.
\vskip0.2cm

Obviously, the convergence of the operators $F_n$ with respect to the continuous dissimilarity index can be obtained also in the more general framework of Theorem \ref{Lp-estimates}. In this case, for any $f \in M(I)$ we immediately get
\be
|1-\text{cSSIM}(f(\cdot),\, F_n(f, \cdot))|\ \leq K\, \left( \tau(f,n^{-1/4})_2\right)^2,
\ee
as $n \to +\infty$, for a suitable $K>0$. 

\section{Other indexes: S-index and PSNR}\label{3.2}

In this section, we provide other indicators for the evaluation of the reconstruction/resizing performances of the algorithm proposed in Section \ref{sec2.33}, as well as for other standard methods, which will be discussed in Section \ref{s4}. These include the \textit{likelihood index}, also known as S-index (see, e.g., \cite{likelihood, likelihood2}), that measures the pointwise deviation between pixel luminances ranging from $0$ to $1$. It is defined as
\begin{equation*}
    \text{S} := \frac{1}{N \cdot M} \sum_{i=1}^{N} \sum_{j=1}^{M}(1-|F_{ij}-G_{ij}|),
\end{equation*}
where $F_{ij}$ and $G_{ij}$ denote the pixel values of two images at position $(i,j)$. A higher value of S, closer to $1$, indicates greater similarity between the images.\\

The \textit{mean square error} (MSE) is another commonly used metric, defined as
\begin{equation}\label{MSE}
\text{MSE} := \frac{1}{N \cdot M} \sum_{i=1}^{N} \sum_{j=1}^{M} |F_{ij} - G_{ij}|^2.
\end{equation}
MSE quantifies the average squared difference between corresponding pixel values of two images, with lower values indicating higher similarity.\\

Furthermore, the \textit{Peak Signal-to-Noise Ratio} (PSNR, see, e.g., \cite{hore2010image}) is defined via the MSE as follows
\begin{equation*}
\text{PSNR} := 20 \cdot \log_{10} \left( \frac{\max{I}}{\sqrt{\text{MSE}}} \right),
\end{equation*}
where $\max{I}$ is the maximum signal measure in the original image. For 8-bit gray scale images $\max{I} = 255$; while for images with pixel values between 0 and 1, $\max{I} = 1$. Such index is commonly used as an indicator to assess the quality of processed images. Indeed, the higher the value of PSNR, the greater the ``similarity'' to the original image, meaning it approaches the original image more closely from a human perceptual standpoint.

\section{Other methods of image rescaling}\label{s4}

Many different interpolation techniques have been developed in the last decades and are used in different settings. To evaluate the performance of the NN algorithm presented in Section \ref{sec2.33}, in the next section we compare it with the commonly used \textit{bilinear} and \textit{bicubic interpolation methods} (see, e.g., \cite{keys1981cubic, smith1981bilinear, liu2013directional, bialecki1994h, getreuer2011linear}). These methods are very classical in digital image processing, and are already implemented in several software and dedicated commands are available in most used programming languages. While they are fast polynomial-based approaches, they have limitations, particularly when handling image edges. 

\begin{remark}
Convergence results related to cSSIM for bilinear and bicubic interpolation methods, specifically providing theoretical bounds on the expected cSSIM as resolution increases, are discussed in \cite{marchetti_santin2022}.
\end{remark}
We will also compare the NN algorithm with u-VPI (\textit{upscaling-de la Vall\'{e}e-Poussin Interpolation}) method based on Lagrange interpolation and filtered de la Vall\'{e}e-Poussin type interpolation, specifically applied at the zeros of Chebyshev polynomials of the first kind. The filter action range provides an additional parameter (denoted by $\theta$, see Section \ref{s5}) that can be suitably regulated to improve the approximation. For the implementation of this method, we will use the openly accessible source code provided by the authors in \cite{occorsio2023image}. Such interpolation technique has also been applied also to three-dimensional images in \cite{occorsio2023filtered}.

\section{Numerical simulations and experiments}\label{s5}

In the numerical experiments, we adopted the following approach. We started by selecting a collection of images with dimensions $N \times M$, which will serve as a reference. Such images are part of the Grayscale Set 1 and 2, provided by the Waterloo Fractal Coding and Analysis Group, and are available in the repository at \url{https://links.uwaterloo.ca/Repository.html}.
\vskip0.2cm

The primary reference set includes four grayscale images: (a) \textit{montage} ($256\times 256$ pixel resolution); (b) \textit{france} ($672\times 496$ pixel resolution); (c) \textit{mountain} ($640\times 480$ pixel resolution); (d) \textit{library} ($464 \times 352$ pixel resolution). The choice of such images is motivated by the fact that we want to make a comparison across images that vary in size, brightness, and texture, while maintaining a common feature: a high variance; particularly, with values exceeding $5e+03$. The latter number came out from the performed numerical experiments. The variance of each one of the considered images are the following: \textit{montage} ($5.0616e+03$), \textit{france} ($5.3166e+03$), \textit{mountain} ($6.4982e+03$), and \textit{library} ($7.9805e+03$).
\vskip0.2cm

We give particular emphasis to the image variances, since as shown in the theoretical Section \ref{3.1} (and confirmed by the numerical evidence), the estimates for the convergence rate are influenced by the constant $c_f$ given in Theorem \ref{th_marchettisantin}. As you can see, such a constant becomes smaller for higher variances, suggesting that the quality of the reconstruction can improve in these situations with respect to the cases characterized by small variances.\\

Then, the selected images were downscaled without interpolation (using the \textit{nearest neighbor method}, \cite{biau2015lectures}) to dimensions $\frac{N}{2} \times \frac{M}{2}$. Subsequently, the downscaled images are upscaled back to their original size using the previously outlined methods. For the algorithm based on NN operators, we considered as $\sigma$ both the logistic and the ramp function, and the following values of $n = 5, 10, 15, 20, 25$ and $30$. Obviously, also others choices of $\sigma$ could be done, but we decide to consider the ones for which we proved in the previous sections the better quantitative estimates.

For the u-VPI method, we set the free parameter $\theta$ equal to $0.5$ following the same approach of the authors in \cite{occorsio2023image}, as we are working in an unsupervised mode where the target image is not available in general.

The results of our experiments are detailed in Tables \ref{Montage}-\ref{Library}, where the best-performing metrics for each image are highlighted in bold for easy identification. To evaluate the PSNR and the (discrete) SSIM with \textsc{Matlab\textsuperscript{\textcopyright}}, we have used its built-in functions \texttt{psnr()} and \texttt{ssim()}, respectively.\\

\begin{table}[tbph]
\centering
    \begin{tabular}{lcccc}
    \hline
    \multicolumn{1}{c}{}      &        & PSNR             & S-index         & SSIM            \\ \hline
    NN ramp function     & $n=10$ & 26.2714          & \textbf{0.9851} & \textbf{0.8004} \\ \hline
    NN logistic function & $n=15$ & \textbf{26.4121} & 0.9843          & 0.7919          \\
                              & $n=30$ & 26.2714          & \textbf{0.9851} & \textbf{0.8004} \\ \hline
    Bilinear interpolation    &        & 24.0623          & 0.9783          & 0.7965          \\ \hline
    Bicubic interpolation     &        & 23.8386          & 0.9773          & 0.7820          \\ \hline
    u-VPI                     &        & 26.3794          & 0.9803          & 0.7680          \\ \hline
    \end{tabular}
    \caption{{\small Comparison of rescaling methods for the image \textit{montage}.}}     \label{Montage}
\end{table}
\begin{table}[tbph]
    \centering
    \begin{tabular}{lcccc}
    \hline
    \multicolumn{1}{c}{}      &        & PSNR             & S-index         & SSIM            \\ \hline
    NN ramp function     & $n=10$ & 19.3282          & \textbf{0.9656} & \textbf{0.7162} \\ \hline
    NN logistic function & $n=25$ & 19.3327          & \textbf{0.9656} & \textbf{0.7162} \\ \hline
    Bilinear interpolation    &        & 18.7431          & 0.9569          & 0.6763          \\ \hline
    Bicubic interpolation     &        & 18.7358          & 0.9570          & 0.6748          \\ \hline
    u-VPI                     &        & \textbf{19.6847} & 0.9624          & 0.6941          \\ \hline
    \end{tabular}
    \caption{{\small Comparison of rescaling methods for the image \textit{france}.} }    \label{France}
\end{table}
\begin{table}[tbph]
\centering
    \begin{tabular}{lcccc}
    \hline
    \multicolumn{1}{c}{}      &        & PSNR             & S-index         & SSIM            \\ \hline
    NN ramp function     & $n=10$ & 16.1870          & \textbf{0.9193} & 0.3940          \\ \hline
    NN logistic function & $n=30$ & 16.1870          & \textbf{0.9193} & 0.3940        \\ \hline
    Bilinear interpolation    &        & 17.3291          & 0.9187          & 0.3782          \\ \hline
    Bicubic interpolation     &        & \textbf{17.3122} & 0.9184          & 0.3750          \\ \hline
    u-VPI                     &        & 16.9448          & 0.9169          & \textbf{0.4025} \\ \hline
    \end{tabular}
    \caption{{\small Comparison of rescaling methods for the image \textit{mountain}.}}     \label{Mountain}
\end{table}
\begin{table}[tbph]
\centering
    \begin{tabular}{lcccc}
    \hline
    \multicolumn{1}{c}{}      &        & PSNR             & S-index         & SSIM            \\ \hline
    NN ramp function     & $n=10$ & 15.4488          & \textbf{0.9306} & \textbf{0.4898} \\ \hline
    NN logistic function & $n=30$ & 15.4488          & \textbf{0.9306} & 0.4896          \\ \hline
    Bilinear interpolation    &        & \textbf{16.6057} & 0.9239          & 0.4668          \\ \hline
    Bicubic interpolation     &        & 16.5855          & 0.9239          & 0.4668          \\ \hline
    u-VPI                     &        & 16.5096          & 0.9242          & 0.4268          \\ \hline
    \end{tabular}
    \caption{{\small Comparison of rescaling methods for the image \textit{library}.}}     \label{Library}
\end{table}	
The experimental results show that the NN algorithm, using both the ramp and the logistic density functions, has a competitive and satisfactory performance. In terms of quality measures, the NN algorithm generally outperforms the other methods considered. 
Additionally, we have also analyzed the CPU time required by each algorithm to process a single image (see, Table \ref{tab_cpu}). On the other hand, the u-VPI method is significantly faster than the others and shows strong performance in terms of PSNR. Moreover, both bilinear and bicubic interpolation methods also exhibit high performance with respect to PSNR.\\
\begin{table}[tbph]
\centering
\begin{tabular}{lccccc}
\hline
\multicolumn{1}{c}{}      &       & \textit{montage}  & \textit{france}   & \textit{mountain} & \textit{library}  \\ \hline
NN ramp function     & $n=5$ & 0.191119          & 0.504026          & 0.468410          & 0.274950          \\
 & $n=10$ & 0.127550 & 0.542623  & 0.505618  & 0.277729 \\
 & $n=15$ & 0.149283 & 0.780755  & 0.606646  & 0.347511 \\
 & $n=20$ & 0.264580 & 1.452352  & 0.677772  & 0.382793 \\
 & $n=25$ & 0.300787 & 1.425558 & 0.796002  & 0.445498 \\
 & $n=30$ & 0.333868 &  1.658628    &   0.975733   & 0.526288 \\ \hline
NN logistic function & $n=5$ & 0.299850          & 1.284863          & 1.232382          & 0.641046          \\
 & $n=10$ & 0.301439 & 1.421931  & 1.385992  & 0.732135 \\
 & $n=15$ & 0.351056 & 1.806180  & 1.682178  & 0.823049 \\
 & $n=20$ & 0.367779 & 2.044200  & 1.866412  & 0.905235 \\
 & $n=25$ & 0.457866 & 2.701355 & 2.141778 & 1.090495 \\
 & $n=30$ & 0.484573 &   2.754693      &   2.468962    & 1.301823 \\ \hline
u-VPI                     &       & \textbf{0.065089} & \textbf{0.018530} & \textbf{0.024777} & \textbf{0.014970} \\ \hline
\end{tabular}
\caption{{\small Comparison of CPU time (sec) for the NN algorithm and the u-VPI method.}}
\label{tab_cpu}
\end{table}

\begin{figure*}[tbph]
    \centering
    \begin{subfigure}{0.48\textwidth}
        \centering
        \includegraphics[width=\linewidth]{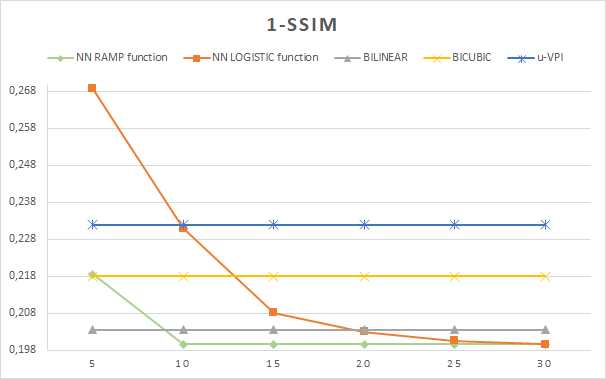}
        \label{graph_montage}
    \end{subfigure}
    \hfill
    \begin{subfigure}{0.48\textwidth}
        \centering
        \includegraphics[width=\linewidth]{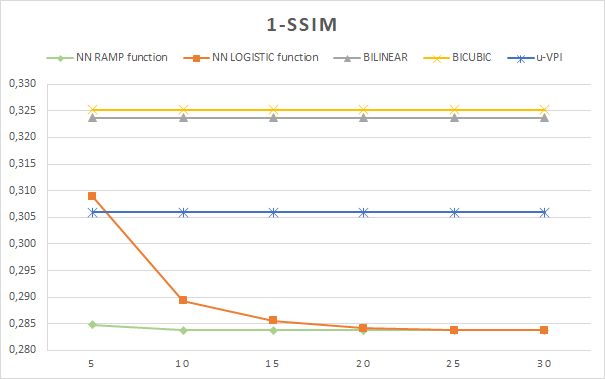}
        \label{graph_france}
    \end{subfigure}
    \caption{{\small Discrete dissimilarity index for the NN algorithm with ramp and logistic functions, plotted with $n = 5, 10, 15, 20, 25, 30$ on the horizontal axis, bilinear and bicubic interpolation, and u-VPI method for the images (a) \textit{montage}, (b) \textit{france}.}}     \label{fig:grafici1}
\end{figure*}

\begin{figure*}[tbph]
    \centering
    \begin{subfigure}{0.48\textwidth}
        \centering
        \includegraphics[width=\linewidth]{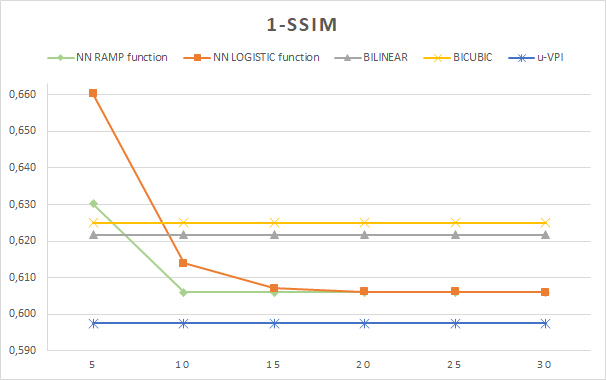}
        \label{graph_montage}
    \end{subfigure}
    \hfill
    \begin{subfigure}{0.48\textwidth}
        \centering
        \includegraphics[width=\linewidth]{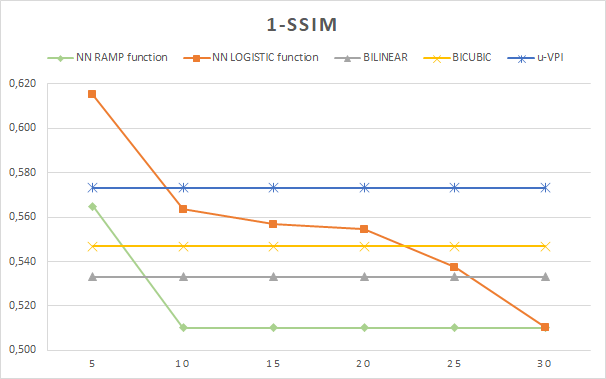}
        \label{graphLibrary.png}
    \end{subfigure}
    \caption{{\small Discrete dissimilarity index for the NN algorithm with ramp and logistic functions, plotted with $n = 5, 10, 15, 20, 25, 30$ on the horizontal axis, bilinear and bicubic interpolation, and u-VPI method for the images (c) \textit{mountain}, (d) \textit{library}.}}     \label{fig:grafici2}
\end{figure*}

The decay trend of the discrete dissimilarity index $1-\text{SSIM}$ is reported in Figure \ref{fig:grafici1} and Figure \ref{fig:grafici2}. The numerical results align with Theorem \ref{corollary_nuovo}. In fact, as $n$ increases, we observe that the dissimilarity index decays when applying the NN algorithm. In particular, the NN algorithm with the ramp function shows a faster decay compared to the one using logistic functions. For completeness, we also present the numerical results for bilinear and bicubic interpolation, as well as the u-VPI methods, even if they do not depend on the parameter $n$.

\begin{remark}
In some preliminary experiments, other than $\sigma_\ell$ and $\sigma_R$ also the hyperbolic tangent  activation function $\sigma_h$ has been considered. However, $\sigma_h$ produced similar results as above, offering no additional advantages. More general sigmoidal functions (see, function $\sigma_\gamma$, with $\gamma>0$ defined on p. 199 of \cite{piconi2024}) could be considered but however they have  a lower convergence rate as showed in Theorem \ref{thnuovo3}.
\end{remark}

To assess robustness and generalization, we extended the tests to six additional greyscale images: \textit{bird} ($256 \times 256$, variance $2.1176e+03$), \textit{boat} ($504 \times 504$, variance $2.6343e+03$), \textit{goldhill} ($256 \times 256$, variance $2.3817 e+03$), \textit{camera} ($256 \times 256$, variance $3.8859e+03$), \textit{slope} ($256 \times 256$, variance $6.8626e+03$), and \textit{mandrill} ($256 \times 256$, variance $1.4778 e+03$). These include both high and low variance cases, allowing for broader evaluation across different image characteristics. For cases where the variance is less than $5e+03$, the performance of the proposed method remains comparable to that of alternative reconstruction techniques, with slightly less good performances in the cases of lower variance.

We also evaluated five colour images from the same repository: \textit{clegg} ($814 \times 880$, variance $6.4331e+03$), \textit{peppers} ($512 \times 512$, variance $4.2271e+3$), \textit{lena} ($512 \times 512$, variance $3.4772e+3$), \textit{frymire} ($1118 \times 1104$, variance $9.8381 e+03$), and \textit{serrano} ($628 \times 794$, variance $6.6334 e+03$).

\begin{figure}[tbph]
    \centering
    \begin{minipage}[t]{0.49\textwidth}
        \centering
        \includegraphics[height=0.25\textheight]{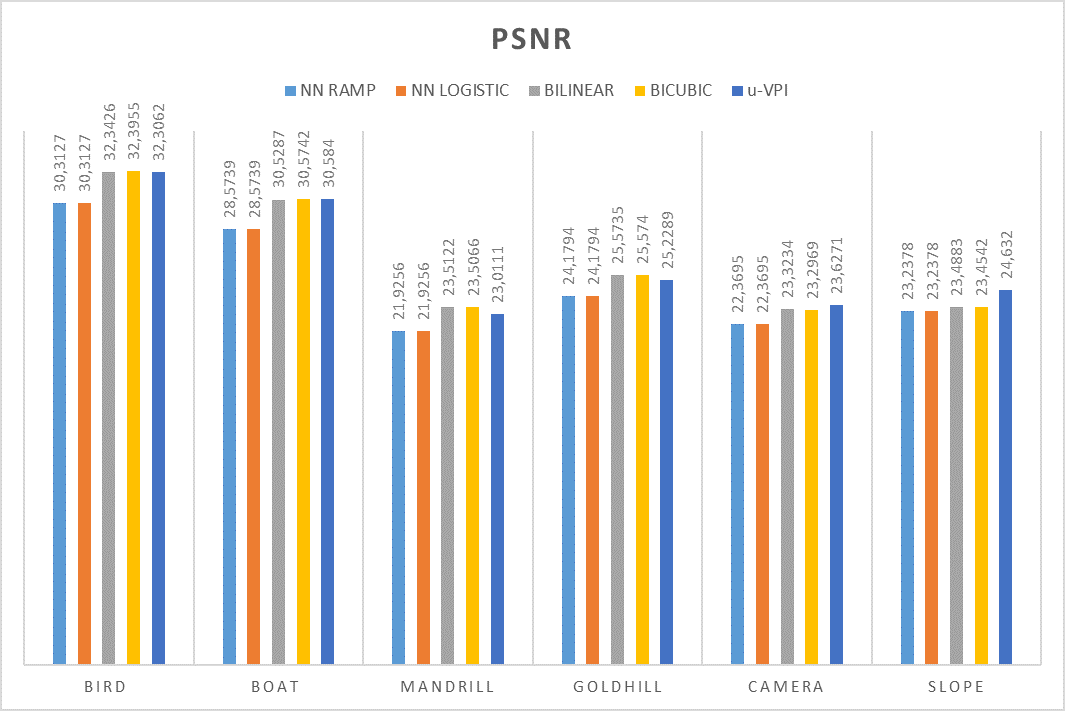}\par\vspace{0.4cm}
        \includegraphics[height=0.25\textheight]{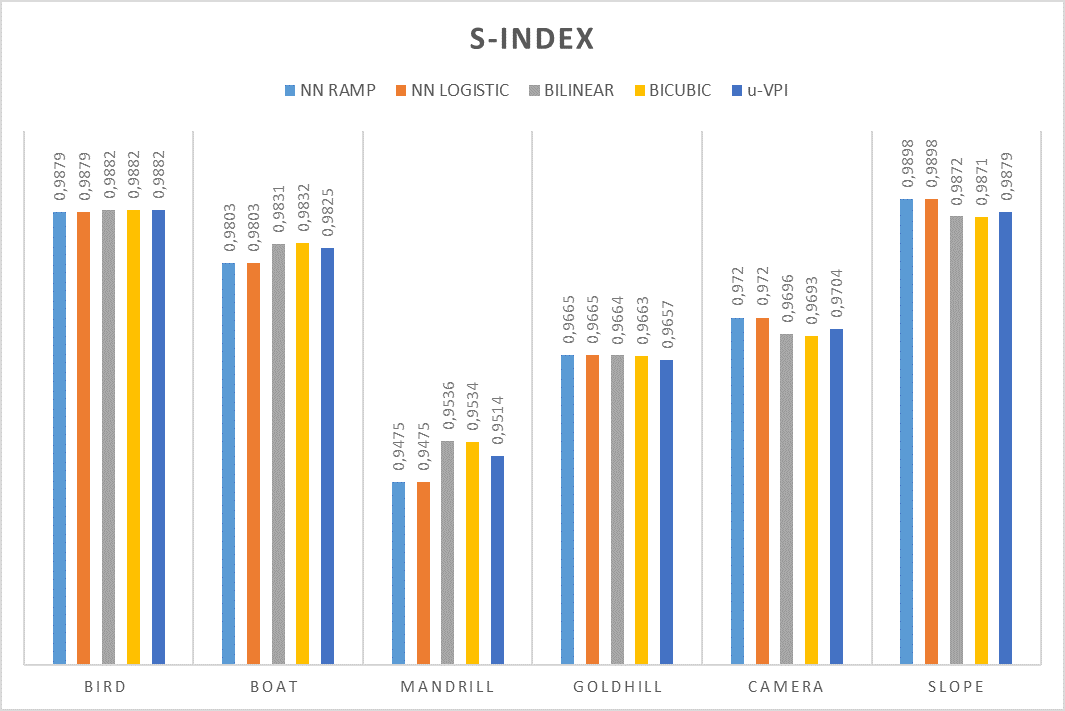}\par\vspace{0.4cm}
        \includegraphics[height=0.25\textheight]{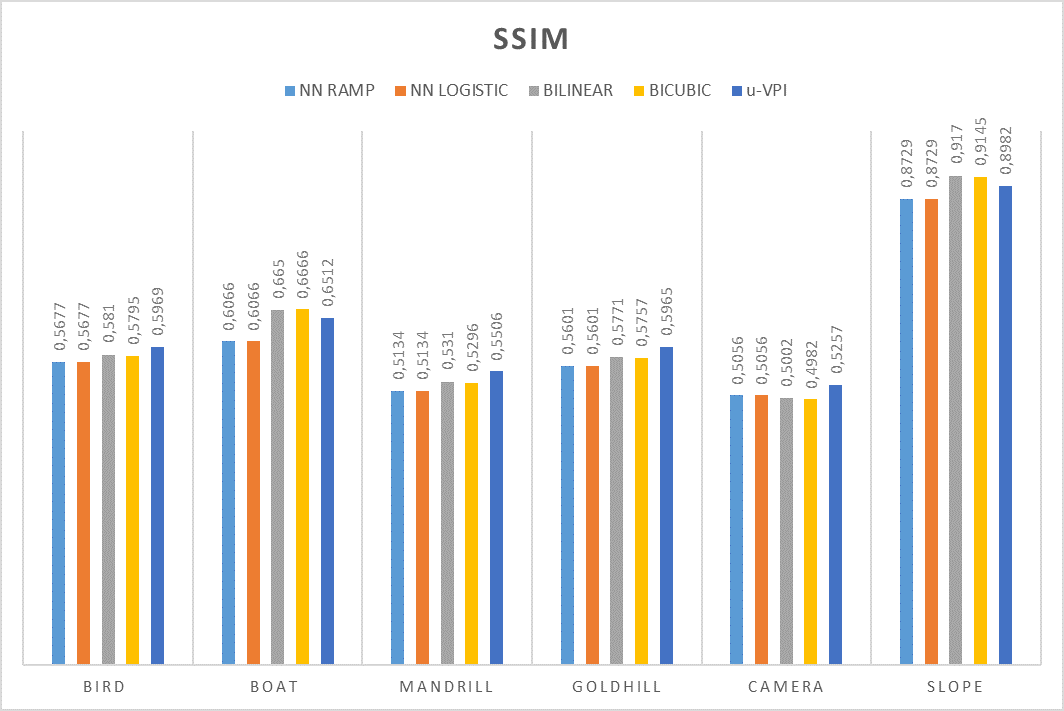}
    \end{minipage}
    \hfill
    \begin{minipage}[t]{0.49\textwidth}
        \centering
        \includegraphics[height=0.25\textheight]{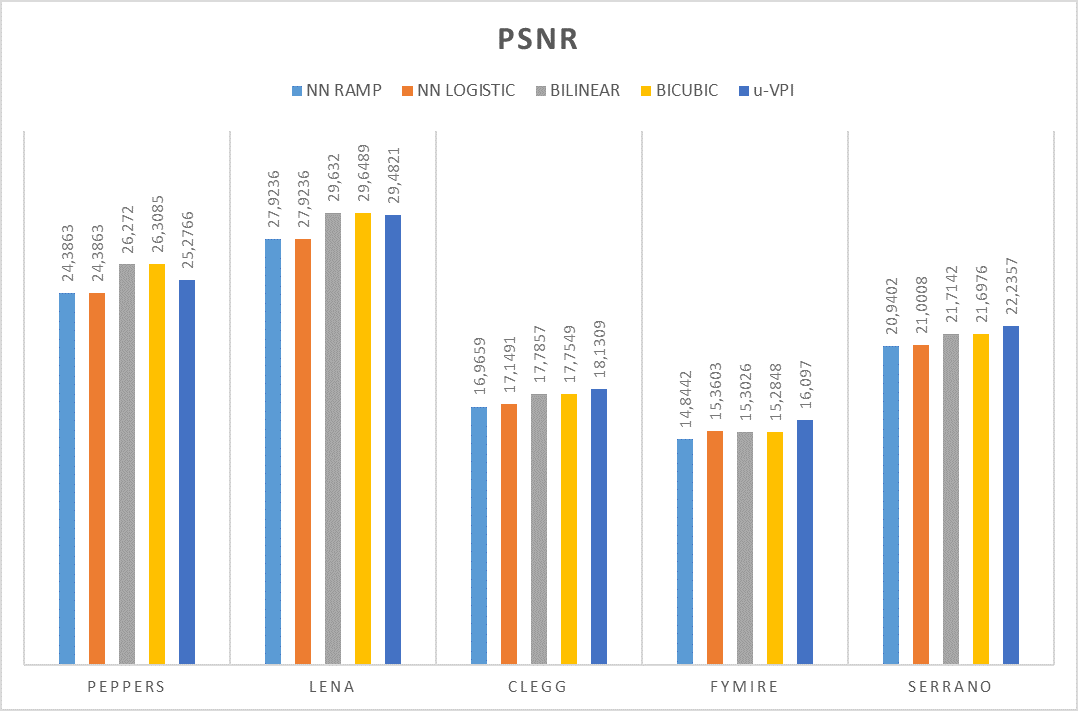}\par\vspace{0.4cm}
        \includegraphics[height=0.25\textheight]{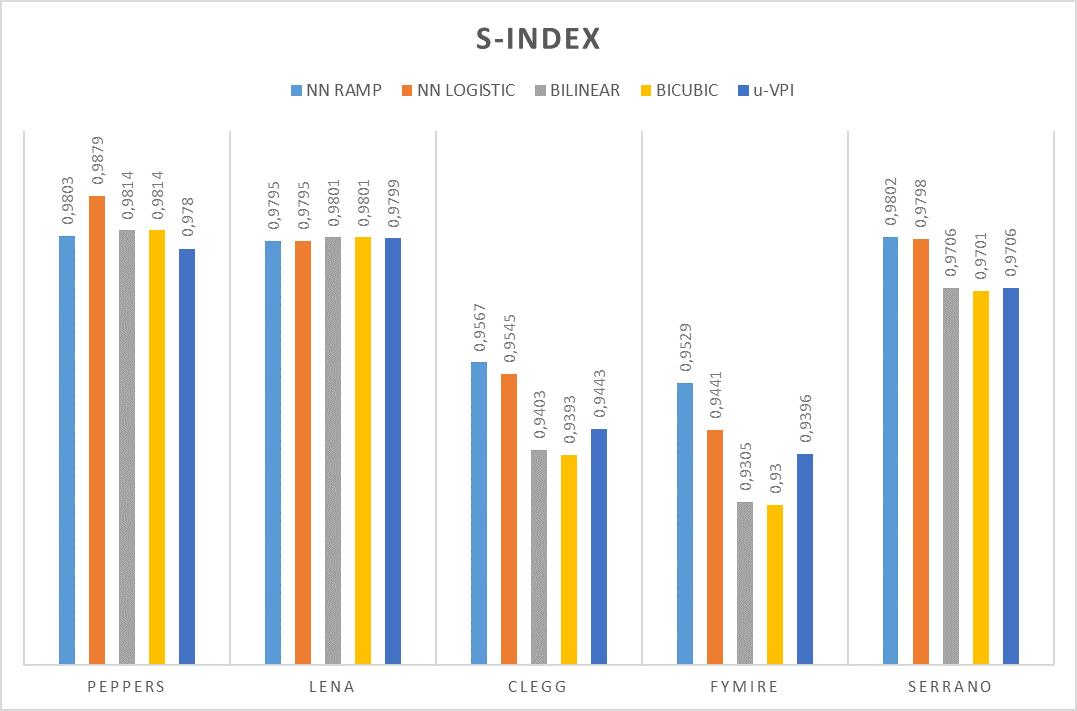}\par\vspace{0.4cm}
        \includegraphics[height=0.25\textheight]{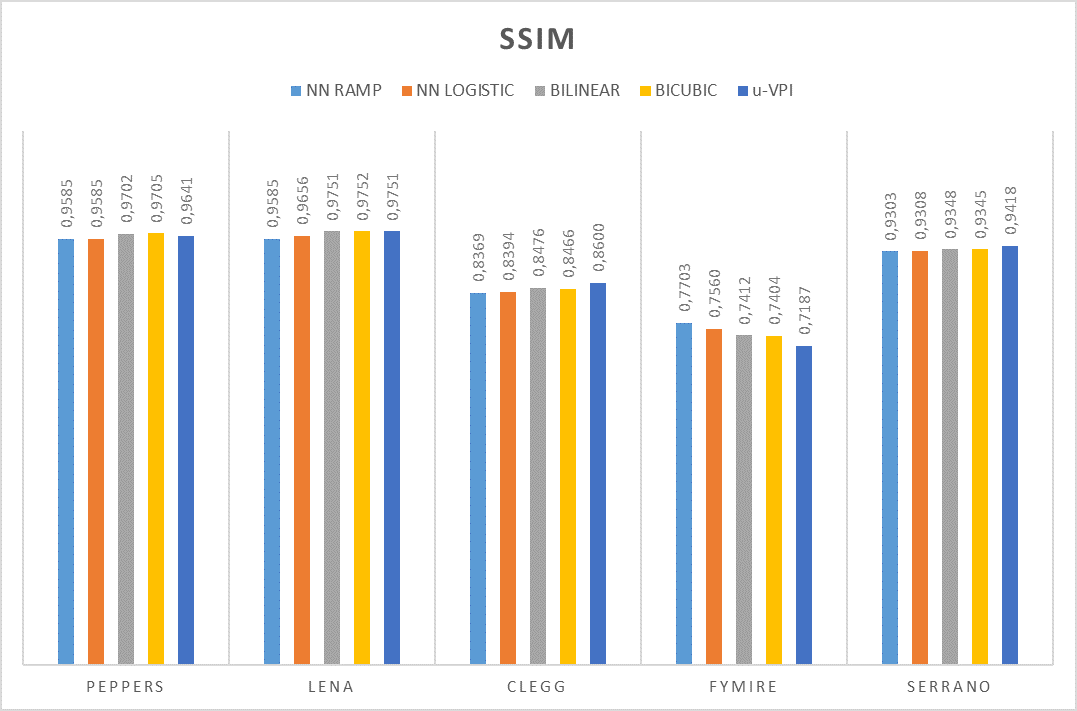}
    \end{minipage}
    \caption{{\small Histogram comparison of the PSNR, S-Index, and SSIM values for different test \textit{grayscale} (left) and \textit{RGB} (right) images (listed along the vertical axis), with the parameter fixed at $n = 30$. Each bar represents the corresponding quality index for a given image, highlighting the performance of the proposed method across multiple metrics.}}     \label{histo_all}
\end{figure}

Due to space limitations, we do not report full tables for all images. Instead, we present aggregated performance metrics in Figure \ref{histo_all}, using histograms of PSNR, SSIM, and S-index values. 

\section{Final remarks, conclusions and future developments}\label{final}

In conclusion, the NN algorithm shows a competitive performance in image rescaling compared to traditional methods such as bilinear, bicubic, and recent u-VPI interpolation method. Specifically, for images with a variance greater than $5e+03$ (where such value has been experimentally determined), the NN approach outperformed the other methods in terms of quality metrics, including SSIM and S-index. While the u-VPI method offers faster processing speed, the NN algorithm provides superior image quality. Experimental data indicate that the NN algorithm performs better for images that cover almost the entire greyscale range (see, Table \ref{valorigrigio} where these data for the considered primary dataset of images have been reported).
\begin{table}[h!]
\centering
\begin{tabular}{cccc}
\hline
\textit{montage} & \textit{france} & \textit{mountain} & \textit{library} \\ \hline
3                & 5               & 146               & 42               \\ \hline
\end{tabular}
\caption{{\small Number of grayscale values not represented by any pixel.}}
\label{valorigrigio}
\end{table}
\vskip0.2cm

\begin{figure*}[tbph]
    \centering
    \begin{subfigure}{0.48\textwidth}
        \centering
        \includegraphics[width=\linewidth]{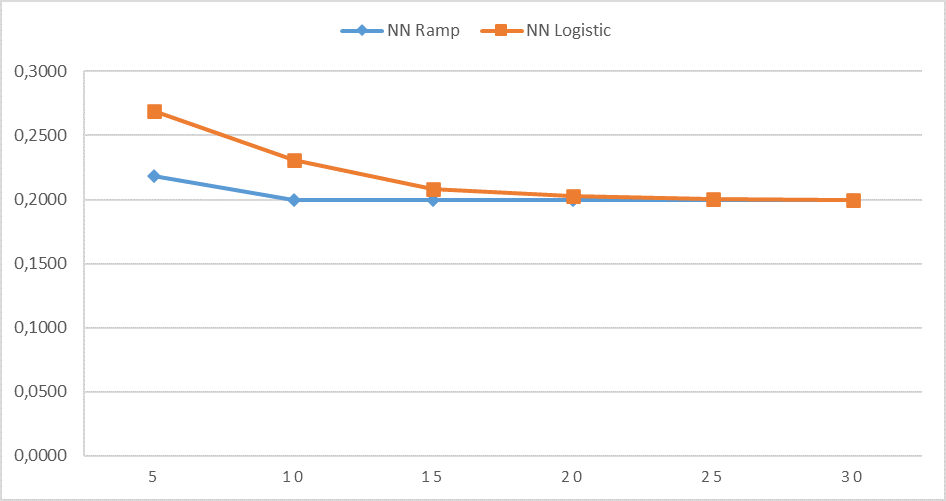}
        \caption{}
        \label{montage_dissimilarity}
    \end{subfigure}
    \hfill
    \begin{subfigure}{0.48\textwidth}
        \centering
        \includegraphics[width=\linewidth]{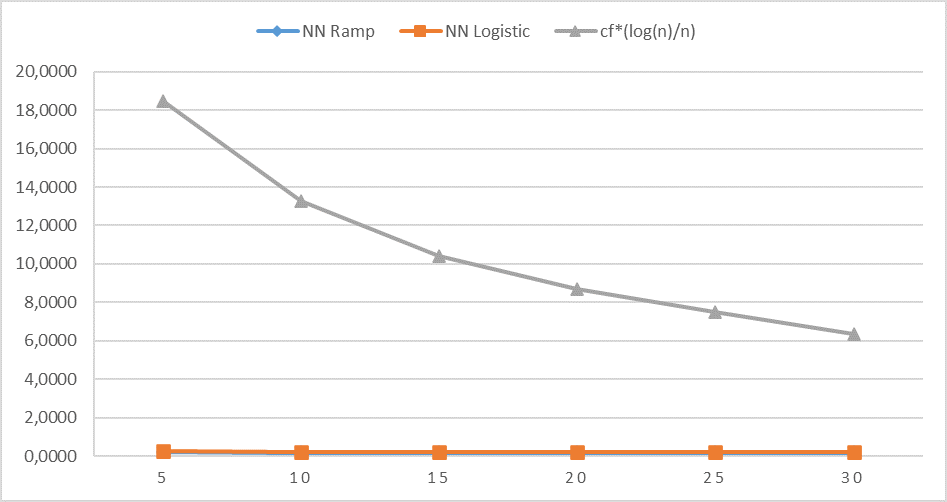}
        \caption{}
        \label{montage_dissimilarity_1}
    \end{subfigure}
    \vspace{0.5cm} 
    \begin{subfigure}{0.48\textwidth}
        \centering
        \includegraphics[width=\linewidth]{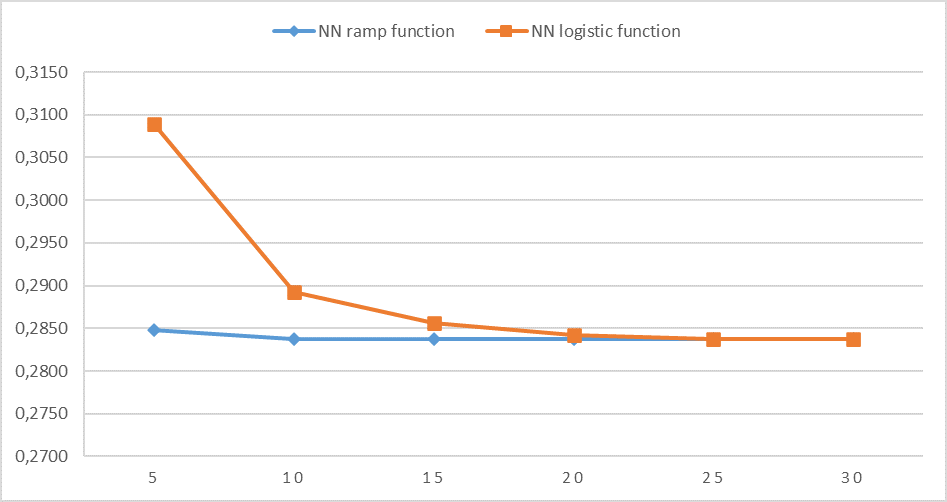}
        \caption{}
        \label{france_dissimilarity}
    \end{subfigure}
    \hfill
    \begin{subfigure}{0.48\textwidth}
        \centering
        \includegraphics[width=\linewidth]{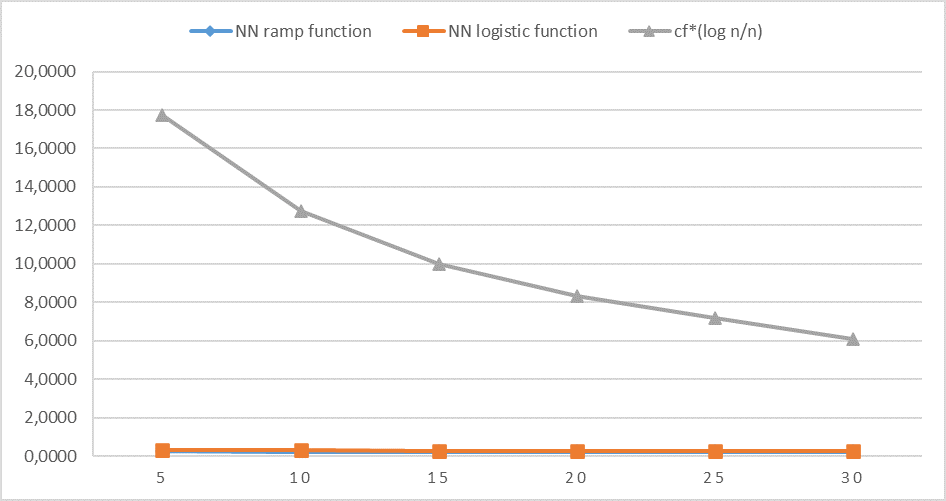}
        \caption{}
        \label{france_dissimilarity_1}
    \end{subfigure}
    \caption{{\small Numerical dissimilarity index $1 - \text{cSSIM}$ obtained by NN algorithm with ramp and logistic functions is shown in (a) for \textit{montage} and in (c) for \textit{france}. In (b) and (d), the same numerical dissimilarity is compared with the theoretical dissimilarity $\frac{cf \cdot \log n}{n}$ as a function of $n = 5, 10, 15, 20, 25, 30$. In both (b) and (d), the numerical dissimilarity remains significantly below the theoretical curve.}}
\label{fig_finale}
\end{figure*}

 Since in the present paper we have devoted a lot of effort to deriving quantitative estimates for the error of approximation, in particular in the special case of image reconstruction, it could be interesting to compare the theoretical analysis with the numerical one. In this direction, we compare the numerical decay of the dissimilarity index with the theoretical asymptotic estimates given in Section \ref{3.1}. According to Corollary~\ref{CorollaryFinal}, the numerical errors seem to decay asymptotically as $\log(n)/n$ (in the cases of $\sigma_\ell$ and $\sigma_R$); such a comparison has been given only for two images (the others behave in the same way): \textit{montage} and \textit{france}, and they can be seen in Figure \ref{fig_finale}. All the parameters involved in generating the plots in Figure \ref{fig_finale} are recorded in Table \ref{tab_parametri}. 
\begin{table}[h!]
\centering
\begin{tabular}{lccccc}
\hline
\multicolumn{1}{c}{} & M   & N   & Variance & Mean     & $c_f$   \\ \hline
\textit{montage}              & 464 & 352 & 5.06e+03 & 100.0395 & 57.8157 \\ \hline
\textit{france}               & 672 & 496 & 5.31e+03 & 99.4246  & 55.4052 \\ \hline
\end{tabular}
\caption{{\small Parameters for the \textit{montage} and \textit{france} images.}}
\label{tab_parametri}
\end{table}
\\
In particular, we assume constants $c_1:=0.03\times 255$ and $c_2:=0.01\times 255$ for the computation of $c_f$, as suggested in the literature (see, e.g., \cite{occorsio2023image}). In particular, in the right plots in Figure \ref{fig_finale} we can visualize the decay of the curves $c_f\, \log(n) /n$, for different values of $n$. These graphs actually show that the estimates established in Corollary~\ref{CorollaryFinal} are (quite large), due to the presence of big constants. Thus, the numerical evidence suggests that the result of Corollary~\ref{CorollaryFinal} could be refined.
\\

Moreover, as already discussed, the value of the constant $c_f$ decreases as the image variance increases (see, e.g., Table \ref{tab_parametri} again). This indicates that, for images with higher variance, the estimates for the dissimilarity index are sharper than those that can be achieved for images with lower variance, in agreement with the numerical observation that the proposed method tends to perform better for high-variance images.\\

Future work should focus on optimizing the NN algorithm by experimenting with different density functions to improve execution speed without compromising its accuracy.

Additionally, we also plan to include deep learning-based super-resolution methods such as EDSR (\textit{Enhanced Deep Residual Networks}), SRGAN (\textit{Super-Resolution Generative Adversarial Network}), and ESRGAN (\textit{Enhanced SRGAN}) in future comparisons, as these models have shown significant improvements in single-image super-resolution  \cite{lim2017enhanced, ledig2017photo, wang2018esrgan}. Given that the NN operator approach is fully deterministic, incorporating such baselines allows a more complete evaluation of its performance within the framework of recent image reconstruction techniques.

At the same time, it is important to highlight that the present NN algorithm works with any single image and does not require training, while methods such as EDSR, SRGAN, etc., require training using large datasets in order to reach high performances. This consideration shows how the approaches are completely different in nature.

\section*{Acknowledgment}
{\small The authors wish to thank the referees for their constructive comments which turn out very useful in order to improve the final version of the present paper.

\noindent  The authors are members of the Gruppo Nazionale per l'Analisi Matematica, la Probabilit\`{a} e le loro Applicazioni (GNAMPA) of the Istituto Nazionale di Alta Matematica (INdAM), of the gruppo UMI (Unione Matematica Italiana) T.A.A. (Teoria dell' Approssimazione e Applicazioni), and of the network RITA (Research ITalian network on Approximation).}

\section*{Funding}
{\small The first and the third authors have been supported within the project PRIN 2022 PNRR: ``RETINA: REmote sensing daTa INversion with multivariate functional modeling for essential climAte variables characterization'', funded by the European Union under the Italian National Recovery and Resilience Plan (NRRP) of NextGenerationEU, under the Italian Minister of University and Research MUR (Project Code: P20229SH29, CUP: J53D23015950001). Moreover, the first and the second authors have been supported within the project PRIN 2022: ``AI- and DIP-Enhanced DAta Augmentation for Remote Sensing of Soil Moisture and Forest Biomass (AIDA)" funded by the European Union under the Italian National Recovery and Resilience Plan (NRRP) of NextGenerationEU, under the MUR (Project Code: 20229FX3B9, CUP: J53D23000660001).}

\section*{Author's contribution} 
{\small All authors contributed equally to this work for writing, reviewing and editing. All authors have read and agreed to the published version of the manuscript.}

\section*{Conflict of interest/Competing interests}
{\small The authors declare that they have no conflict of interest and competing interest.}

\section*{Copyright}
{\small All the images are contained in the repository in \url{https://links.uwaterloo.ca/Repository.html} \cite{Repository}, and they belong to the Grayscale Set 1 and 2, and to the Colour Set (The Waterloo Fractal Coding and Analysis Group).\\ The source code for the proposed u-VPI method from \cite{occorsio2023image} is publicly accessible at \url{https://github.com/ImgScaling/VPIscaling}. Permission is granted to use, copy, or modify the documentation for educational and research purposes without fee. }

\section*{Availability of data and material and Code availability}
{\small All the data generated for this study were stored in our laboratory and are not publicly available. Researchers who wish to access the data directly contacted the corresponding author.}


\end{document}